\documentclass[12pt,a4paper,reqno]{amsart}
\usepackage{amssymb}
\usepackage{latexsym}
\usepackage{exscale}

\numberwithin{equation}{section}
\newtheorem{theor}{Theorem}[section]
\newtheorem{lemma}[theor]{Lemma}
\newtheorem{defi}[theor]{Definition}
\newtheorem{corol}[theor]{Corollary}
\newtheorem{remark}[theor]{Remark}
\newtheorem{prop}[theor]{Proposition}
\newtheorem{defin}[theor]{Definition}

\newcommand{\re}{\mathbb{R}}
\newcommand{\co}{\mathbb{C}}

\newcommand{\RR}{\mathbb{R}}
\newcommand{\NN}{\mathbb{N}}
\newcommand{\B}{\mathbb{B}}

\newcommand{\J}{\mathcal{J}}
\newcommand{\K}{\mathcal{K}}

\newcommand{\W}{\mathcal{W}}

\newcommand{\T}{\mathcal{T}}
\newcommand{\D}{\mathcal{D}}
\newcommand{\C}{\mathcal{C}}

\newcommand{\call}{\mathcal{L}}

\newcommand{\X}{\mathcal{X}}

\newcommand{\ep}{\varepsilon}

\newcommand{\bigchi}{\mathop{\mathchoice%
{\mbox{\Large$\chi$}}{\mbox{\large$\chi$}}{\mbox{\normalsize$\chi$}}%
{\mbox{\small$\chi$}}}\nolimits}

\renewcommand{\emptyset}{\mbox{\rm \O}}
\renewcommand{\Re}{{\rm Re}\,}

\def\div{\mathop{\rm div}}
\def\Int{\mathop{\rm Int}}
\def\supp{\mathop{\rm supp}}

\newcommand{\aver}[1]{-\hskip-0.46cm\int_{#1}}

\newcommand{\offx}[2]{\mathcal{O}\big(L^{#1}-L^{#2}\big)}

\newcommand{\off}[2]{\mathcal{O}\big(L^{#1}(\mu)-L^{#2}(\mu)\big)}
\newcommand{\offw}[2]{\mathcal{O}\big(L^{#1}(w)-L^{#2}(w)\big)}
\newcommand{\full}[2]{\mathcal{F}\big(L^{#1}(\mu)-L^{#2}(\mu)\big)}

\newcommand{\fullx}[2]{\mathcal{F}\big(L^{#1}-L^{#2}\big)}

\newcommand{\dec}[1]{\Upsilon\!\left(#1\right)}

\newcommand{\expt}[1]{e^{\textstyle #1}}

\begin{document}
\allowdisplaybreaks

\title[Weighted norm inequalities and elliptic
operators]{Weighted norm inequalities,  off-diagonal estimates and
elliptic operators
\\[.2cm]
{\footnotesize Part II:  Off-diagonal estimates \\ on spaces of
homogeneous type}}

\author{Pascal Auscher}

\address{Pascal Auscher
\\
Universit\'e de Paris-Sud et CNRS UMR 8628
\\
91405 Orsay Cedex, France} \email{pascal.auscher@math.p-sud.fr}

\author{Jos\'e Mar{\'\i}a Martell}

\address{Jos\'e Mar{\'\i}a Martell
\\
Instituto de Matem\'aticas y F{\'\i}sica Fundamental
\\
Consejo Superior de Investigaciones Cient{\'\i}ficas
\\
C/ Serrano 123
\\
28006 Madrid, Spain}

\address{\null\vskip-.7cm and\vskip-.7cm\null}

\address{Departamento de Matem\'aticas \\ Universidad Aut\'onoma de Madrid \\
28049 Madrid, Spain } \email{chema.martell@uam.es}

\thanks{This work was partially supported by the European Union
(IHP Network ``Harmonic Analysis and Related Problems'' 2002-2006,
Contract HPRN-CT-2001-00273-HARP). The second author was also
supported by MEC ``Programa Ram\'on y Cajal, 2005'' and by MEC Grant
MTM2004-00678.}

\date{\today}


\begin{abstract}
This is the second part of a series of four articles on weighted
norm inequalities, off-diagonal estimates and elliptic operators. We
consider a substitute to the notion of pointwise bounds for kernels
of operators which   usually is a measure of decay. This substitute
is that of off-diagonal estimates expressed in terms of local and
scale invariant $L^p-L^q$ estimates. We propose a  definition in
spaces of homogeneous type that is stable under composition.  It is
particularly well suited to semigroups. We  study the case of
semigroups generated by elliptic operators.
\end{abstract}

\maketitle

\null\vskip-2cm\null

\begin{quote}
{\scriptsize\tableofcontents}
\end{quote}

\section{Introduction}

Since its discovery by Gaffney \cite{Gaf} for the heat equation on
Riemannian manifolds, one can show that most semigroups  generated by
elliptic operators satisfy the so-called $L^2$ off-diagonal estimates.
They are  of the form
\begin{equation}\label{eq:offL2L2}
\|T_t (\bigchi_{E} f)\|_{L^2(F)} \le C \expt{-\frac{c\,
d^2(E,F)}{t}} \|f\|_{L^2(E)},
\end{equation}
valid for all $t>0$ and all $f\in L^2$ whenever $E, F$ are closed
sets and $d(E,F)$ an appropriate distance on sets. This estimate is
relevant when $t $ is smaller than a time $t_0$ comparable to
$d^2(E,F)=d^2$:  before  $t_0$, the heat, which we imagine with a
Gaussian distribution,  has not had  enough  time to diffuse from $E$
to  give a  significant contribution on  $F$, hence the decay in
$t/d^2$ which explains the terminology ``off-diagonal.'' These
estimates are instrumental in many applications of semigroups. For
example, they are the main technical tool (for the resolvent instead
of the semigroup) in the proof of the Kato conjecture \cite{AHLMcT}.
Whether one can improve integrability properties depends on other
arguments. Such an improvement,
 called hypercontractivity (see, e.g.  \cite{Da3}), is usually linked to some kind of Sobolev
embedding.

A strong form of  off-diagonal estimates  is the Gaussian upper
bound, that is,  a pointwise control of the kernel  of $T_{t}$ by
Gaussians:
\begin{equation}
\label{eq:gub}
|K_{t}(x,y)|\le Ct^{-n/2} \expt{-\frac{c\, d^2(x,y)}{t}}.
\end{equation}
This behavior appears frequently and has yielded in the 1990's a
number of  beautiful results  on  independence of sectors of
analyticity for semigroups, independence of $L^p$-spectrum  as $p$
varies for their generators, maximal regularity problems, \dots. An
account on all this as well as a documented bibliographical list can
be found in a recent survey by  Arendt \cite{Arendt}.

The power of $t$ in front of the Gaussian factor appears in
homogeneous situations where the volume of  balls is comparable to a
power  of their radii. In this case, such an estimate implies
$L^1-L^\infty$ boundedness of $T_{t}$ known as the ultracontractivity
property. There are geometric situations, such as  Riemannian
manifolds or weighted measures, where the volume of a ball is not
comparable to a power of its radius. In this case, the Gaussian upper
bound becomes
\begin{equation}
\label{eq:gub1}
|K_{t}(x,y)|\le \frac C {\sqrt{{\rm Vol}\, (B(x,\sqrt t))\, {\rm Vol}\, (B(y,\sqrt t))}\,}\, \expt{-\frac {c\,  d^2(x,y)}t},
\end{equation}
 and this no longer implies   ultracontractivity.  This estimate has to be treated as some sort of local and scale invariant $L^1-L^\infty$ bound which can be called $L^1-L^\infty$ off-diagonal estimates; it is still the improvement of regularity in the scale of Lebesgue spaces that matters, even if it is local.

The Gaussian upper bound \eqref{eq:gub} is equivalent to
\begin{equation}\label{eq:offL1Linfty}
\|T_t (\bigchi_{E} f)\|_{L^\infty(F)} \le C t^{-n/2}
\expt{-\frac{c\,  d^2(E,F)}{t}} \|f\|_{L^1(E)}
\end{equation}
for all $t>0$, $f\in L^1$ and all closed sets $E, F$.  Interpolation
between  \eqref{eq:offL2L2} and \eqref{eq:offL1Linfty} yield
intermediate  $L^p-L^{p'}$ conditions of the same type. Hence, one
can formulate a definition    for arbitrary  $p,q$ with $1\le p<q\le
\infty$ which we call here $L^p-L^q$ full off-diagonal estimates.
Such conditions  appear naturally  in absence of Gaussian upper
bounds \eqref{eq:gub}.  Davies showed that such a generalization
already leads to improved results on independence of sectors of
analyticity and independence of $L^p$-spectrum  as $p$ varies
\cite{Da4}.

In another direction, we owe to Blunck and Kunstmann the fundamental observation that this
notion of  $L^p-L^q$ off-diagonal estimates for $p<q$  permits to
develop a theory of  singular ``non-integral'' operators for which
one can formulate $L^p$ boundedness criteria  for $p$ in arbitrary
intervals in the absence of information on the kernels (pointwise
bounds and even bounds in mean). Let us mention their weak type
$(p,p)$ criterion for  $1<p<2$ in absence of  kernels and assuming
weak type $(2,2)$, similar to the generalization by Duong and
McIntosh \cite{DMc} of H\"ormander's result \cite{Hor} when $p=1$ in
presence of kernel bounds. See their series of papers
\cite{BK1,BK2,BK3}. See also \cite{HM} for related ideas. One can
find in Fefferman's work \cite{Fef}  the essence of such a criterion
but no explicit statement was needed because it was in a situation
where one can split operators in pieces with localized smooth
kernels. In the same spirit, \cite{ACDH} proposes a strong type
$(p,p)$  criterion for one $p>2$ (not all) via good-$\lambda$
inequalities (we also refer the reader  to the first article of the
series \cite{AM1} where we generalize this to weighted norm estimates). In
\cite{Aus}, all these ideas are presented in the Euclidean setting
and applied to some singular  ``non-integral'' operators arising from
elliptic operators. This yields optimal ranges of exponents $p$ for
$L^p$ boundedness;  the weighted norm extension for this application
is the purpose of \cite{AM3}. However, we observe that the  $L^p-L^q$
full  off-diagonal estimates when $p<q$ used by Blunck and Kunstmann,
even if they are stable under composition (which is a natural
property when working with semigroups),  is somewhat irrealistic in
general situations for at least three reasons: they implies global
$L^p-L^q$ boundedness of $T_{t}$, they do not imply $L^p-L^p$
boundedness of $T_{t}$, and they do not  pass to weighted estimates.
So the notion of $L^p-L^q$ off-diagonal estimates, generalizing
\eqref{eq:gub1}, on a space  of homogeneous type  for one-parameter
families of operators  needs a definition.

Such a  definition should  only involve
balls and annuli  and make clear that they are two parameters involved,  the radius of balls  and the  parameter of the family,  linked by a scaling rule independently on the  location of the balls.
Some examples suggest possible definitions (called here strong or
mild off-diagonal estimates) but they are no longer stable under composition in a
general context.  Hence, the price to pay for  stability is a somewhat
weak definition (in the sense that we can not be greedy in our
demands). Nevertheless, it covers examples of the literature on
semigroups.  Furthermore,   in spaces of homogeneous type with
polynomial volume growth (that is, the measure of a ball is
comparable to  a power of its radius, uniformly over centers and
radii) it  coincides with all other definitions.  This is also the case for more general volume growth conditions, such as the one for some Lie groups with a local dimension and a dimension at infinity. Eventually,  it
is  operational for proving weighted estimates  in \cite{AM3},
which was the main motivation for developing this material. Since
it is of independent interest, we present  it  here in a
separate article, which can be read independently of the other papers of our series.

In Section \ref{sec:ode}, we introduce our definition of $L^p-L^q$ off-diagonal
estimates on balls for one-parameter (such as time) families of
operators and  state the main properties:  for $p=1, q=\infty$, it is
equivalent to the  Gaussian upper bound \eqref{eq:gub1}  for the
associated kernels, and for arbitrary $p,q$, it implies $L^p$
boundedness of each operator,  is stable under composition, and
passes to the weighted case  (proofs are given later in Section~\ref{sec:proofs}).
As mentioned,  we  discuss  in Section \ref{sec:other} other ``expected'' stronger
definitions for off-diagonal estimates. The proof that all our
definitions coincide in spaces with polynomial  volume growth  is in
Section \ref{sec:proofs}. We then present in Section \ref{sec:propandsemi} the application to
semigroups: we establish a correspondence between   the interval of
exponents $p$ for $L^p$-boundedness of $T_{t}$ and the set of
exponents $(p,q)$ for which   $L^p-L^q$ off-diagonal estimates on
balls hold  when the latter set is not empty. This correspondence
remains  true for sectorial analytic extension of semigroups with
independence of angles. Eventually,  we show how unweighted
off-diagonal estimates imply weighted ones   for appropriate
$A_{\infty}$ weights.  In Section \ref{sec:casestudy}, we apply all this and  describe
unweighted and weighted off-diagonal estimates on balls of semigroups
$\{e^{-t\, L}\}_{t>0}$ and their gradient $\{\sqrt t\, \nabla \,
e^{-t\, L}\}_{t>0}$ for a class of elliptic operators $L$ in $\RR^n$.

\section{Off-diagonal estimates on balls}\label{sec:ode}

\subsection{Setting and notation}\label{sec:setting}

Let $(\X,d,\mu)$ be a space of homogeneous type,  which is a (non empty) set $\X$
endowed with a distance $d$ (it could even be a quasi-distance but we restrict to this situation for simplicity)  and a non-negative Borel measure $\mu$ on
$\X$ such that the doubling condition
\begin{equation}\label{doubling}
\mu(B(x,2\,r)) \le C_0\, \mu(B(x,r))<\infty,
\end{equation}
holds for all $x\in\X$ and $r>0$, where $B(x,r)=\{y\in\X:d(x,y)<r\}$.

Throughout this paper we use the following notation: for every
ball $B$, $x_B$ and $r_B$ are respectively its center and its radius,
that is, $B=B(x_B,r_B)$. Given $\lambda>0$, we will write
$\lambda\,B$ for the $\lambda$-dilated ball, which is the ball with
the same center as $B$ and with radius $r_{\lambda\,B}=\lambda\,r_B$.

If $C_0$ is the smallest constant for which the measure $\mu$
verifies the doubling condition $\eqref{doubling}$, then $D=\log_2
C_0$ is called the doubling order of $\mu$ and we have that
$\mu(\lambda\,B)\le C_\mu\,\lambda^D\,\mu(B)$, for every ball $B$ and
for every $\lambda\ge 1$.

Given a ball $B$ we set $C_j(B)=2^{j+1}\,B\setminus
2^j\,B$ for $j\ge 2$; $C_1(B)=4\,B$ and also
$\widehat{C}_1(B)=4\,B\setminus 2\,B$. We set
$$
\aver{B} h\,d\mu
=
\frac1{\mu(B)}\,\int_B h\,d\mu,
\qquad\qquad
\aver{B^c} h\,d\mu
=
\frac1{\mu(2\,B)}\,\int_{B^c} h\,d\mu,
$$
and for $j\ge 1$
$$
\aver{C_j(B)} h\,d\mu
=
\frac1{\mu(2^{j+1}\,B)}\,\int_{C_{j}(B)} h\,d\mu.
$$
The last notation  can be seen as the average on $2^{j+1}B$ of
$\bigchi_{C_{j}(B)} h$, where we  denote by $\bigchi_E$ the indicator
function of a set $E$. It is not necessarily the case that $2^{j+1}\,B$ and
$C_{j}(B)$ have comparable masses with constant independent of $B$ and $j$ (for example, when $C_{j}(B)=\emptyset$) so it is safer to
divide out by the mass of the larger set (which is never 0 unless $\mu=0$) and fortunately, this is the quantity arising in computations. Although this is not needed in this work, let us mention a reasonable sufficient condition insuring this comparability (see \cite{AM1}  for a proof): Assume that there exists  $\ep\in (0,1)$ such that for any ball $B\subset \X$,  $(2-\ep)\, B \setminus B \ne \emptyset$. Then, $ \mu(2\,B) \lesssim \mu(2\,B\setminus B)$ for any ball $B$, where the implicit constants are independent of $B$.

For shortness we write
$\dec{s}=\max\{s,s^{-1}\}$ for $s>0$. We use the symbol $A \lesssim
B$ for $A \le C B$ for some constant $C$ whose value is not important
and independent of the parameters at stake.

\subsection{Definition and comments}

\begin{defi}\label{defi:off-d:weights}
Given $1\le p\le q\le \infty$, we say that a family $\{T_t\}_{t>0}$ of sublinear
operators  satisfies $L^{p}(\mu)-L^{q}(\mu)$
off-diagonal estimates on balls, which  by an abuse of notation will
be denoted $T_t \in\off{p}{q}$, if there exist $\theta_1, \theta_2>0$
and $c>0$  such that for every $t>0$ and for any ball $B$, setting $r=r_{B}$,
\begin{equation}\label{w:off:B-B}
\Big(\aver{B} |T_t( \bigchi_B \, f) |^{q}\,d\mu\Big)^{\frac 1 q}
\lesssim
\dec{\frac{r}{\sqrt{t}}}^{\theta_2} \,\Big(\aver{B}
|f|^{p}\,d\mu\Big)^{\frac 1 p };
\end{equation}
and, for all $j\ge 2$,
\begin{equation}\label{w:off:C-B}
\Big(\aver{B}|T_t( \bigchi_{C_j(B)}\, f) |^{q}\,d\mu\Big)^{\frac 1 q}
\lesssim
2^{j\,\theta_1} \dec{\frac{2^j\,r}{\sqrt{t}}}^{\theta_2}\,
\expt{-\frac{c\,4^{j}\,r^2}{t}} \,
\Big(\aver{C_j(B)}|f|^{p}\,d\mu\Big)^{\frac 1 p }
\end{equation}
and
\begin{equation}\label{w:off:B-C}
\Big(\aver{C_j(B)}|T_t( \bigchi_B \, f) |^{q}\,d\mu\Big)^{\frac 1 q}
\lesssim
2^{j\,\theta_1} \dec{\frac{2^j\,r}{\sqrt{t}}}^{\theta_2}\,
\expt{-\frac{c\,4^{j}\,r^2}{t}}
\,\Big(\aver{B}|f|^{p}\,d\mu\Big)^{\frac 1 p }.
\end{equation}
\end{defi}

\noindent\textit{Comments.}
\begin{list}{$\theenumi.$}{\usecounter{enumi}\leftmargin=1cm
\labelwidth=0.7cm\itemsep=0.3cm\topsep=.1cm}

\item When $q=\infty$ one has to change the $L^q$-norms by
the corresponding essential suprema.

\item  $T_{t}$ may only be defined on a subspace of $L^p(\mu)$  provided this subspace
is stable under truncation by indicator functions  of measurable sets (balls would suffice for the definition but  measurable sets is needed for interpolation). In this case, it is understood that
the definition applies to functions $f$ in this subspace.

\item Even though our definition makes sense when $p\ge q\ge 1$, we restrict ourselves to  $p\le q$ to stress  the
regularizing effect in the scale of Lebesgue spaces.

\item H\"older's inequality implies $
\off{p}{q}\subset\off{p_1}{q_1}
$
 for all $p_1,q_1$ with $p\le p_1\le q_1\le
q$.

\item  $T_t\in \off{p}{q}$ with $p<q$ does not imply that $T_{t}$ is bounded from $L^p(\mu)$ into $L^q(\mu)$.

\item If $T_{t}$ is linear and defined on a dense subspace of
$L^p(\mu)$, then $T_{t}\in \off{p}{q}$ if and only if $T_{t}^*\in
\off{q'}{p'}$ where $T_{t}^*$ is the dual operator for the duality
form $\int_{\X} f g \, d\mu$.

\item Given two Banach spaces $B_1$ and $B_2$, this
definition, with the corresponding changes, is also valid for
operators taking $B_1$-valued functions into $B_2$-valued functions.

\item The chosen ``time-space'' scaling $\frac r {\sqrt t}$ is irrelevant. It can be changed
at will to $\frac r {t^\gamma}$ for any $\gamma>0$ simply by changing $T_{t}$ to $T_{t^{2\gamma}}.$ This scaling corresponds to semigroups of second order operators, which is our main application in \cite{AM3}.

\item The profile $s\mapsto e^{-c\,s^2}$ can be replaced by
any non increasing  $g:\re^+\to  \re^+$
 such that for all $\theta\ge 0$, we have
$s^\theta\,g(s)\rightarrow 0$ as $s\rightarrow \infty$. For
example $g(s)=e^{-c\,s^\alpha}$ with $c, \alpha>0$ is acceptable. The value of $c$ has no interest to us provided it remains non negative. Thus, we will freely use the same $c$ from line to line. Profiles with sufficiently large polynomial decay $(1+s)^{-p}$
works as well, but $p$   would have to be
adjusted to each application.

\item For $s\le 1$,  $\dec{s}^{\theta_{2}} e^{-cs^2}$ is comparable to  $s^{-\theta_{2}}$. Since $s$ is to be replaced by
$\frac {2^j\, r}{\sqrt t}$, it is
curious at first sight that we allow such negative powers;  imposing  positive powers for small $s$ seems more natural. We were forced into negative powers to obtain  stability  under
composition. See Lemma \ref{lemma:sum} below for the technical reason.
   Fortunately, this apparently weak behavior is
sufficient for our applications in
\cite{AM3}.

\item One can replace $\theta_1$  by $\theta_1+\alpha$ for any
$\alpha\ge 0$ and the same happens with $\theta_2$. In fact,
making $\theta_{2}\ge \theta_{1}$, one obtains an equivalent
definition by  replacing $2^{j\,\theta_1}
\dec{\frac{2^j\,r}{\sqrt{t}}}^{\theta_2}\expt{-\frac{c\,4^{j}\,r^2}{t}}$
by  an expression of the form
$\dec{\frac{\,r}{\sqrt{t}}}^{\theta_2}\expt{-\frac{c\,4^{j}\,r^2}{t}}$
(up to changing the  $c$'s). We stick to the first formulation for
simplicity in some calculations but this is a first indication
that the value of the exponent $\theta_{1}$ is irrelevant.

\item Definition \ref{defi:off-d:weights} is given in terms of dyadic
annuli but an equivalent definition can be written in terms of
$a$-adic annuli for all $a>1$. See the proof of Lemma \ref{lemma:1-ann} for
a possible argument.
\end{list}

\subsection{The case $p=1$ and $q=\infty$}

We first state that the definition for $p=1$ and $q=\infty$
coincides with the usual pointwise Gaussian decay of the
introduction.

\begin{prop}\label{prop:1infty} Assume that the operators $T_{t}$, $t>0$,
are linear.
Then $T_t \in\off{1}{\infty}$ if and only if there exist constants
$C,c>0$ and for each $t>0$, a measurable function $K_{t}$ on
$\X\times \X$ such that $T_{t}f(x)=\int_{\X}K_{t}(x,y)f(y)\,
d\mu(y)$ holds for  almost every $x\in \X$ whenever $f\in
L^1(\mu)$ and for almost  every $(x,y)\in \X\times \X$,
\begin{equation}
\label{eq:gub2}
|K_{t}(x,y)| \le \frac C {\mu(B(x,\sqrt t))}\, \expt{-\frac {c\,  d^2(x,y)}t}.
\end{equation}
\end{prop}

The (easy) proof in Section \ref{sec:proof:prop:1infty} shows the role of the independence of $r$ and $t$ in the condition $T_t \in\off{1}{\infty}$ and of the scaling rule $r/\sqrt t$.  It also shows the irrelevance of the exponents $\theta_{1},\theta_{2}$ in this case.  Note that    the doubling condition on $\mu$ implies that  the Gaussian expressions in \eqref{eq:gub1} and \eqref{eq:gub2} are comparable  up to changing the constants $C,c>0$.

\subsection{Uniform boundedness and stability under composition}

We state here the most important features of this notion:
Off-diagonal estimates on balls imply uniform boundedness  and are
stable under composition.

\begin{theor}\label{theor:uniform-compo}

\

\begin{list}{$(\theenumi)$}{\usecounter{enumi}\leftmargin=1cm
\labelwidth=.7cm\labelsep=0.25cm \itemsep=0.25cm \topsep=.3cm
\renewcommand{\theenumi}{\alph{enumi}}}

\item If $T_t\in \off{p}{p}$ then $T_t:L^{p}(\mu)\longrightarrow
L^{p}(\mu)$ is bounded uniformly on $t$.

\item Let $1\le p\le q \le r\le \infty$. If $T_t\in\off{q}{r}$ and
$S_t\in\off{p}{q}$ then $T_t\circ S_t\in \off{p}{r}$. Furthermore if
$\theta_1$, $\theta_2$ are the exponents appearing in Definition
\ref{defi:off-d:weights} for $T_t$ and $\gamma_1$, $\gamma_2$ are the
ones for $S_t$, then the composition $T_t\circ S_t \in\off{p}{r}$
satisfies the corresponding inequalities with
$\alpha_1=\max\{\theta_1,\gamma_1,D/r\}$ and
$\alpha_2=\max\{\theta_1,\theta_2\}+\max\{\gamma_1+D/q,\gamma_2\}$.
\footnote[2]{When $\theta_1\ne \theta_2$ and $\gamma_1+D/q\ne
\gamma_2$, the value of $\alpha_{2}$ is correct. Otherwise,
$\alpha_{2}$ is any number strictly bigger than this value,  see
Remark \ref{remark:powers-diff} below.}

\end{list}

\end{theor}

In $(b)$, if $T_{t}$ and $S_{t}$ are defined on subspaces, this
result is understood in the sense that one restricts to functions $f$
for which $T_{t}\circ S_{t} f$ is well-defined.

The  proof of this result can be found in Sections
\ref{sec:proof:theor:uniform-compo-a} and
\ref{sec:proof:theor:uniform-compo-b}.

\subsection{Weighted off-diagonal estimates}\label{subsec:weighted-off}

Weighted off-diagonal estimates on balls
with weights in the Muckenhoupt class $A_\infty$ can be obtained from the
off-diagonal estimates on balls with respect to the underlying
measure.

We use the following notation: given a weight $w$  we consider the
measure $dw=w\,d\mu$, so $w(E)=\int_E dw=\int_E
w\,d\mu$ and also $L^p(w)=L^p(w\,d\mu)$. Notice
that the notation for the averages used before depends on the measure
we are using and so
$$
\aver{B} h\,d\mu
=
\frac1{\mu(B)}\,\int_B h\,d\mu,
\qquad\qquad
\aver{B} h\,dw
=
\frac1{w(B)}\,\int_B h\,dw
=
\frac1{w(B)}\,\int_B h\,w\,d\mu,
$$
the same happens for ``averages''  on $B^c$ and $C_j(B)$ as defined in Section \ref{sec:setting}.

Let $w\in A_\infty$ (we recall some basic facts about $A_{p}$ and $RH_{s}$
weights in Appendix \ref{App:Weights}). Since $(\X,d,\mu)$ is a space of homogeneous
type, the measure $w$ is
doubling and  $(\X,d,w)$ is also a space of homogeneous type. Hence, off-diagonal estimates make sense in that space.

\begin{prop}\label{prop:off-unw-w}
Let $1\le p_0<q_0\le \infty$ and $T_t\in\off{p}{q}$ for all $p$, $q$
with $p_0<p\le q< q_0$. Then, for all $p$, $q$ with $p_0<p\le q< q_0$
and for any $w\in A_{\frac{p}{p_0}}\cap  RH_{(\frac {q_{0}}q)'}$ we
have that $T_t\in \offw{p}{q}$.
\end{prop}

The proof of this result can be found in Section
\ref{sec:proof:prop:off-unw-w}.

\section{Other types of off-diagonal estimates}\label{sec:other}

\subsection{Full off-diagonal estimates}

In the case where $(\X,d,\mu)$ is the usual Euclidean space with Lebesgue measure
or more generally, a group with polynomial volume growth (we say that  $(\X,d,\mu)$ has  polynomial volume growth when $\mu(B(x,r)) \sim r^n$ for some $n>0$ and uniformly for all $x\in \X$ and $r>0$), one encounters more precise off-diagonal estimates. This yields a possible definition in spaces of homogeneous type.

\begin{defin}\label{def:full} Let $(\X,d,\mu)$ be a space of homogeneous type.
Let $1\le p\le q \le \infty$. We say that a family $\{T_t\}_{t>0}$ of
sublinear operators satisfies  $L^p(\mu)-L^q(\mu)$  \textbf{full
off-diagonal estimates}, in short $T_{t}\in \full{p}{q}$,  if for
some constant $ \theta\ge0$, with $\theta\ne 0$ when $p<q$, for all
closed sets $E$ and $F$, all $f$ and all $t>0$ we have
\begin{equation}\label{eq:offLpLq-w}
\Big(\int_{F}|T_t (\bigchi_{E}\,  f)|^q\, d\mu\Big)^{\frac 1 q} \lesssim  t^{-\theta}
\expt{-\frac{c d^2(E,F)}{t}} \Big( \int_{E}|f|^p\, d\mu\Big)^{\frac 1 p }.
\end{equation}
\end{defin}

Again, the operators are defined on some subspace $\D$ that is stable
under truncation by indicators of measurable sets. Full off-diagonal
estimates appear when dealing with  semigroups of second order
elliptic operators (see \cite{Gaf, Da2, LSV, Aus} \dots). The most
studied case   is when $p=1$ and $q=\infty$ which means that the
kernel of $T_{t}$ has pointwise Gaussian upper bounds (see \cite{Aro,
FS,  Cou, VSC, Da3, Rob, AMcT, AT, ArTe, DER}\dots). If one considers
higher order operators, then $t$ changes to some  positive power of
$t$ and the Gaussian to other exponential like function (if $t$
denotes time) (see \cite{Da4, AT}).  Our $t$ here may not be the
usual time scale and the Gaussian may be changed also. We stick to
this case to keep the presentation simple.

When $\X=\RR^n$, the usual value of $\theta$ (given our choice of
``space-time'' scaling) is  $\theta=\frac 1 2
\big( \frac n p - \frac n q\big)$. See the proof of Proposition \ref{prop:full-off}.

Here is a list of simple and known facts whose proofs will be left to
the reader.
\begin{list}{$(\theenumi)$}{\usecounter{enumi}\leftmargin=.8cm
\labelwidth=0.7cm\itemsep=0.1cm\topsep=.2cm
\renewcommand{\theenumi}{\alph{enumi}}}
\item $T_{t}\in \full{p}{q}$  implies $T_{t}$ bounded from $L^p(\mu)$ to $L^q(\mu)$.

\item If $p\le r \le q$, $S_{t}\in \full{p}{r}$ and $T_{t}\in \full{r}{q}$, then  $T_{t}\circ S_{t} \in \full{p}{q}$.
\end{list}

Let us compare this definition with the previous one. On the one
hand, if $q=p$ then $ \full{p}{p}$ easily implies  $\off{p}{p}$. The
converse is true and follows from the proof of Proposition
\ref{prop:full-off}, $(b)$, stated below. No further condition on the
space  is needed at this point.

On the other hand, these notions cease to be comparable when $p<q$ without further information on the space $\X$.
Assume that  $T_{t}\in \full{p}{q}$.  If   $E=F=B$ is a ball, then
$$
\Big(\aver{B}|T_t (\bigchi_{B}\,  f)|^q\, d\mu\Big)^{\frac 1 q} \lesssim  t^{-\theta} \, \mu(B)^{\frac 1 p -\frac 1 q}
\Big( \aver{B}|f|^p\, d\mu\Big)^{\frac 1 p }.
$$
Unless there is some control from above of $\mu(B)$ by  a power of
the radius of $B$, we cannot conclude that $T_{t}\in \off{p}{q}$.
Similarly unless there is such a control, one cannot conclude that
$T_{t}$ uniformly bounded on $L^r(\mu)$ for $p\le r \le q$.
Eventually, if $p \le p_{1}< q_{1}\le q$, we do not know if $T_{t}\in
\full{p_{1}}{q_{1}}$.

As  $L^p(\mu)-L^q(\mu)$ full off-diagonal estimates when $p<q$  imply
$L^p(\mu)-L^q(\mu)$ boundedness but not $L^p(\mu)$ boundedness, this
is not an encountered notion  on a general space of homogenous type.
For example, the heat semigroup $e^{-t\Delta}$ on functions  for
general Riemannian manifolds with the doubling property is not $L^p-L^q$ bounded when $p<q$ unless, as Proposition
\ref{prop:full-off} will show,    the measure of any ball is bounded
below by a power of  its radius.

Here is a statement that connects both notions. The proof  is postponed to Section
\ref{sec:proof:prop:full-off}.

\begin{prop}\label{prop:full-off}
Let $(\X,d,\mu)$ be a space of homogeneous type and $1\le p<q\le \infty$.
\begin{list}{$(\theenumi)$}{\usecounter{enumi}\leftmargin=.8cm
\labelwidth=0.7cm\itemsep=0.2cm\topsep=.2cm
\renewcommand{\theenumi}{\alph{enumi}}}

\item Assume that $\X$ has volume growth at most polynomial, that is,
$\mu(B(x,r)) \lesssim r^n$ for some $n>0$ and uniformly for all $x\in
\X$ and $r>0$. If  $T_{t}\in \full{p}{q}$ with   exponent $\theta$ in
\eqref{eq:offLpLq-w}  equal to $\frac 1 2 \big( \frac n p - \frac n
q\big)$  then $T_{t}\in \off{p}{q}$.

\item Assume that $\X$ has  volume growth at least   polynomial, that
is, $ \mu(B(x,r)) \gtrsim r^n$ for some $n>0$ and uniformly for all
$x\in \X$ and $r>0$. If $T_{t}\in \off{p}{q}$ then $T_{t}\in
\full{p}{q}$ with   exponent $\theta$ in \eqref{eq:offLpLq-w}  equal
to $\frac 1 2 \big( \frac n p - \frac n q\big)$.

\end{list}
 \end{prop}

Let us go a little further. We say that a  space of homogeneous
type is of  $\varphi$-growth if $\mu(B(x,r)) \sim \varphi(r)$
uniformly for $x\in \X$ and $r>0$, where $\varphi$ is a
non-decreasing function on $(0,\infty)$. Remark that the fact that
space is of homogeneous type implies that $\varphi$ is doubling in
the sense that $\sup_{{r>0}} \frac{\varphi(2r)}{\varphi(r)}
<\infty$.\footnote[2]{We  think that the discussion can be
extended somehow to spaces with exponential growth, but this is
beyond the scope of the present article.}  A particular important
example is  the Heisenberg group equipped with Riemannian distance
and Haar measure: in this case,  $\varphi(r)\sim r^d$ for $r\le 1$
and $\varphi(r)\sim r^D$ for $r\ge 1$, the exponents $d>0$ and
$D>0$ being called  respectively  its local dimension  and its
dimension at infinity. Call ``$L^p(\mu)- L^q(\mu)$ full
off-diagonal estimates of type $\varphi$''  the estimates of
Definition  \ref{def:full} with $t^{-\theta}$ replaced by
$\varphi\big(\sqrt t\, \big) ^{\frac 1 q - \frac 1 p }$ for all
$t>0$.

\begin{prop}\label{prop:full-off:2D}
Let $(\X,d,\mu)$ be a  space of homogeneous type of $\varphi$-growth and $1\le p< q\le \infty$.
 Then $T_{t}\in \off{p}{q}$ if and only if $T_{t}$ satisfies $L^p(\mu)- L^q(\mu)$ full off-diagonal estimates of type $\varphi$.
 \end{prop}

In other words, the scaling $\frac r {\sqrt t}$ contained in the off-diagonal estimates on balls plus the volume growth completely rule the function of $\sqrt t$ in the full off-diagonal estimates.

The proof of this result is postponed to Section
\ref{sec:proof:prop:full-off}.

\

Our last remark is that full off-diagonal estimates  do not pass to
weighted measures as well: for example, in $\RR^n$, the power weights
$w(x)= |x|^{- \alpha}$ for $0<\alpha<n$ are  neither  with polynomial
growth from below or above nor with $\varphi$-growth.

\subsection{Mild off-diagonal estimates on balls}

As in Section \ref{subsec:weighted-off}  full off-diagonal estimates imply some (but not full)  off-diagonal weighted estimates for an appropriate class of weights.  Assume that $\X=\RR^n$ equipped with Lebesgue measure. With the same arguments (see Section \ref{sec:proof:prop:off-unw-w})  we
obtain that if  $1\le p_0<q_0\le \infty$ and $T_t$ satisfies
$L^p(dx)-L^q(dx)$ full off-diagonal estimates \eqref{eq:offLpLq-w}  for all $p$, $q$ with $p_0<p\le
q< q_0$ and $\theta=\frac 1 2 \big( \frac n p - \frac n q\big)$, then, for all $p$, $q$ with $p_0<p\le q< q_0$ and for any
$w\in A_{\frac{p}{p_0}}\cap  RH_{(\frac {q_{0}}q)'}$ we have that
\begin{equation}\label{eq:weightedoff}
\Big(\aver{B'} \bigchi_{K'}\,
|T_{t}f|^q\, dw \Big)^{\frac1q}
\lesssim
c(t,B,B',K,K')
\Big(\aver{B}
\bigchi_K\,|f|^{p}\, dw \Big)^{\frac1p}
\end{equation}
where
$$
c(t,B,B',K,K') = \left( \frac {\sqrt t}{r_{B'}}\right)^{\frac n{q_{1}}}
\left( \frac  {r_{B}}{\sqrt t}\right)^{\frac n{p_{1}}}
 \expt{-\frac{c\, d^2(K',K)}{t}}
$$
whenever $B,B'$ are balls, $K,K'$ are respective compact subsets, $f$
bounded with support in $K$, $t>0$ and $p_{1}, q_{1}$ are some numbers
chosen with $p_{0}<p_{1}<q_{1}<q_{0}$. For $q=\infty$, the left hand
side of \eqref{eq:weightedoff} is understood as the essential
supremum on $B'$. If
we specialize to the three cases of Definition
\ref{defi:off-d:weights}, namely, 1) $B=B'=K=K'$, 2) $B=K,
B'=2^{j+1}B, K'=C_{j}(B)$ and 3) the symmetric case of 2), we obtain
\eqref{w:off:B-B}, \eqref{w:off:C-B}, \eqref{w:off:B-C}  with
$\theta_{1}= n/q_{1}$ and $\theta_{2}= n/p_{1}-n/q_{1}> 0$ and
$\dec s$ is replaced by $s$.

This leads us to another definition of off-diagonal estimates in a general context.

\begin{defin}
Let $(\X,d,\mu)$ be a space of homogeneous type. Let $1\le p\le q \le
\infty$. We say that a family $\{T_t\}_{t>0}$ of sublinear operators
satisfies  $L^p(\mu)-L^q(\mu)$  \textbf{mild off-diagonal estimates
on balls}  if there exist real numbers $ \theta_1,  \theta_2 \ge 0$,
$c>0$ with $\theta_{2}>0$ when $p<q$    such that \eqref{w:off:B-B},
\eqref{w:off:C-B} and \eqref{w:off:B-C} hold with $s$ replacing
$\dec s$.
\end{defin}

\begin{remark}\rm
In replacing $\dec s$ by $s$ then one cannot enlarge $\theta_{2}$ at
will as in the definition of off-diagonal estimates on balls. Hence,
the restriction that $\theta_{2}$ should be non negative when $p<q$
seems meaningful.
\end{remark}

This is clearly  stronger  than Definition
\ref{defi:off-d:weights} since we impose the power of $s$ to be
positive even for small $s$ (see comment 10 after Definition
\ref{defi:off-d:weights}). However, stability under composition is
unclear. If $S_{t}$ satisfies $L^p(\mu)-L^q(\mu)$ mild
off-diagonal estimates on balls and $T_{t}$ satisfies
$L^q(\mu)-L^r(\mu)$ mild off-diagonal estimates on balls, then we
do not know whether $T_{t}\circ S_{t}$ satisfies
$L^p(\mu)-L^r(\mu)$ mild off-diagonal estimates on balls. Of
course, $T_{t}\circ S_{t} \in \off{p}{r}$ (hence,  under
$\varphi$-growth there is stability)

We may have lost too much information in passing from full
off-diagonal estimates  to mild off-diagonal estimates on balls,
hence the lack of stability.  In particular, we restricted attention
to balls while the closed sets $E$ and $F$ in \eqref{eq:offLpLq-w} could be unbounded.

\subsection{Strong off-diagonal estimates on balls}
The
following result will suggest an even stronger definition.

\begin{prop}\label{prop:full-strong}
Assume that $(\X,d,\mu)$ is the usual Euclidean space $\RR^n$ with
Lebesgue measure.  Fix $1 \le p_{0} < q_{0}\le \infty$. Assume that
$\{T_t\}_{t>0}$ satisfies $L^p(dx)-L^q(dx)$ full off-diagonal estimates  for
all $p,q$ with $p_{0}<p\le q < q_{0}$ and $\theta=\frac 1 2 \big( \frac n p - \frac n q\big)$.  Fix $p,q$ with $p_{0}<p\le q
< q_{0}$ and assume that $w \in A_{\frac p {p_{0}}} \cap RH_{(\frac
{q_{0}}q)'}$. Let  $B$ be a ball and set $r=r_{B}$. Then for all $f$,
\begin{equation}
\label{eq:full-strong:B-2Bc}
\Big( \aver{(2B)^c} |T_{t}(\bigchi_{B}\, f)|^q\, dw \Big)^{\frac 1 q} \lesssim
 \Big( \frac  {r}{\sqrt t}\Big)^{\beta}\expt{-\frac{c\, r^2 }t}  \Big
(\aver{B} |f|^{p} \,\, dw \Big)^{\frac 1 p }
\end{equation}
and
\begin{equation}
\label{eq:full-strong:2Bc-B}
\Big (\aver{B} |T_{t}(\bigchi_{(2B)^c} \, f )|^q\, dw \Big)^{\frac 1 q} \lesssim
 \Big( \frac  {r}{\sqrt t}\Big)^{\gamma}\expt{-\frac{c \,r^2 }t}  \Big
(\aver{(2B)^c} |f|^{p} \, dw \Big)^{\frac 1 p }
\end{equation}
with $\beta, \gamma\ge 0$  and non zero when $p<q$.
\end{prop}

The proof of this result is postponed until Section
\ref{sec:proof:prop:full-strong}.

\begin{defin}
Let $(\X,d,\mu)$ be a space of homogeneous type. Let $1\le p\le q \le \infty$. We say that a family $\{T_t\}_{t>0}$ of sublinear operators
 satisfies  $L^p(\mu)-L^q(\mu)$  \textbf{strong off-diagonal estimates on
balls}  if there exist real numbers $ \alpha\ge 0$ with $\alpha>0$ when $p<q$ and $c>0$ such that  for any ball $B$ and any $t>0$, setting $r=r_{B}$, and any $f$ in an appropriate space $\D$, \begin{equation}\label{w:st-off:B-B}
\Big(\aver{B} |T_t( \bigchi_B \, f) |^{q}\,d\mu\Big)^{\frac 1 q}
\lesssim
\Big({\frac{r}{\sqrt{t}}}\Big)^{\alpha}  \,\Big(\aver{B}
|f|^{p}\,d\mu\Big)^{\frac 1 p };
\end{equation}
\begin{equation}\label{w:st-off:B-2Bc}
\Big(\aver{(2B)^c} |T_t( \bigchi_B\, f) |^{q}\,d\mu\Big)^{\frac 1 q}
\lesssim
\Big({\frac{r}{\sqrt{t}}}\Big)^{\alpha} \expt{-\frac{c\,r^2}{t}}\,\Big(\aver{B}
|f|^{p}\,d\mu\Big)^{\frac 1 p };
\end{equation}
\begin{equation}\label{w:st-off:2Bc-B}
\Big(\aver{B} |T_t( \bigchi_{(2B)^c} \, f) |^{q}\,d\mu\Big)^{\frac 1 q}
\lesssim
\Big({\frac{r}{\sqrt{t}}}\Big)^{\alpha} \expt{-\frac{c\,r^2}{t}} \,\Big(\aver{(2B)^c}
|f|^{p}\,d\mu\Big)^{\frac 1 p }.
\end{equation}
\end{defin}

It is clear that strong off-diagonal estimates on balls imply mild
off-diagonal estimates on balls: for instance, to get the analog of
\eqref{w:off:C-B}, we write $\widetilde{B}=2^{j-1}\,B$ and note that
$C_j(B)\subset (2\,\widetilde{B})^c$. So we apply
\eqref{w:st-off:2Bc-B} with $\widetilde{B}$ and then we obtain
\eqref{w:off:C-B} with $s$ replacing $\dec{s}$, $\theta_2=\alpha$ and
$\theta_1=D/q$. The same can be done in the other cases.

It is also interesting to compare the last two inequalities of this
definition with the ones in Lemma \ref{lemma:Bc}: again $\dec s$ is
replaced by $s$.  We also stress that such a definition implies the
partial $L^p(\mu)-L^q(\mu)$ boundedness inequalities
$$
\Big(\aver{B} |T_t f|^{q}\,d\mu\Big)^{\frac 1 q}
\lesssim
\Big({\frac{r}{\sqrt{t}}}\Big)^{\alpha} \,\Big(\frac 1{\mu(B)} \int_{\X}
|f|^{p}\,d\mu\Big)^{\frac 1 p }$$
and
$$
\Big(\frac 1 {\mu(B)}\int_{\X} |T_t( \bigchi_B \, f) |^{q}\,d\mu\Big)^{\frac 1 q}
\lesssim
\Big({\frac{r}{\sqrt{t}}}\Big)^{\alpha} \,\Big(\aver{B}
|f|^{p}\,d\mu\Big)^{\frac 1 p }.
$$
Had we put an estimate from $B^c$ to $B^c$ similar to
\eqref{w:st-off:B-B} then we would derive global $L^p(\mu)-L^q(\mu)$
boundedness, which is not realistic when $p<q$ in a non polynomial growth situation.

For this reason precisely, strong off-diagonal estimates do not
compose well. The best we can say (even if we allow the exponent
$\alpha$ to take different values in \eqref{w:st-off:B-B},
\eqref{w:st-off:B-2Bc}, \eqref{w:st-off:2Bc-B})  is: If $S_{t}$
satisfies $L^p(\mu)-L^q(\mu)$ strong off-diagonal estimates on balls
and $T_{t}$ satisfies $L^q(\mu)-L^r(\mu)$ strong off-diagonal
estimates on balls then $T_{t}\circ S_{t}$ satisfies  $L^p(\mu)-
L^r(\mu)$ mild off-diagonal estimates on balls (this can be obtained
easily following the proof of $(b)$ in Theorem
\ref{theor:uniform-compo} and using the definition of the strong
off-diagonal estimates on balls in place of Lemma \ref{lemma:Bc}). Again, assuming $\varphi$-growth, there is stability under composition, passing via full off-diagonal estimates.

\

In conclusion,  using only balls, complements of balls and annuli for defining off-diagonal estimates (instead of closed sets) forces us into an apparently weak definition to have
stability under composition. But under a polynomial (or $\varphi$-) growth, all these notions are the same.

\section{Propagation and semigroups}\label{sec:propandsemi}

We are interested in values of $p,q$ for which $L^p-L^q$
off-diagonal estimates on balls hold, especially when there is a
regularizing effect, that is, when $p<q$.

\subsection{Propagation property}\label{sec:propagation}

Let $\T=\{T_{t}\}_{t>0}$ be a family of sublinear operators defined
on a space $\D$ contained in all $L^p(\mu)$ that is stable under
truncation by indicator functions of measurable sets.

Let  $\widetilde \J(\T)$ be the interval of all exponents $p\in [1,\infty]$ such that
$T_{t}$ is bounded uniformly with respect  $t$ on  $L^p(\mu)$.

We introduce the set
$$
{\mathcal{O}}(\T)=\{ (p,q)\in [1,\infty]^2\, ; \,  p<q,\ T_t\in\off{p}{q}\}.
$$
If we set $\C(\T)=\{(\frac 1 p ,\frac 1 q)\, ;\, (p,q) \in {\mathcal{O}}(\T)\}$, then by interpolation, it is a convex set contained in $\{(u,v)\in [0,1]^2\, ; \, u> v\}$.

The relation between ${\mathcal{O}}(\T)$ and  $\widetilde \J(\T)$ is the following.
If $(p,q) \in {\mathcal{O}}(\T)$ then the interval $[p,q]$ is contained in $\widetilde \J(\T)$.
This fact is a consequence of Theorem \ref{theor:uniform-compo}, part $(a)$.

Also, if $ {\mathcal{O}}(\T)\ne \emptyset $,  then for
$p\in\Int\widetilde \J(\T)$,\footnote[2]{If $E$ is a subset of
$[1,\infty]$ with lower and upper bound $p,q$ then we set $\Int
E=(p,q)=\{t\in \RR\, ; \, p<t<q\}$, which is the interior of $E\cap
\RR$ in $\RR$.} there exists  $q=q(p)>p$ such that $T_{t}\in
\off{p}{q}$. In other words, $L^p(\mu)$ boundedness improves into
some off-diagonal estimates on balls with increase of exponent or,
differently, off-diagonal estimates on balls for one pair $(p,q)$
propagate to pairs $(p,q(p))$ for  all $p\in \Int\widetilde
\J(\T)$.\footnote[3]{Here, we see $p$ as the exponent in the source
space. It could also be taken as the exponent of the target space:
for $q\in\Int\widetilde \J(\T)$, there exists  $p=p(q)<q$ such that
$T_{t}\in \off{p}{q}$.} Indeed, let $p\in \Int\widetilde \J(\T) $.
Let $(q,r)\in {\mathcal{O}}(\T)$. If $p=q$ we have finished.
Otherwise, we have that $p,q\in \widetilde \J(\T)$,  and since $p$ is
in the interior, there exists $\tilde p\in \widetilde \J(\T)$ such
that  $p$ lies in the open interval between $\tilde p$ and $q$. The
$L^{\tilde p}(\mu)$ boundedness implies that $T_{t}$ satisfies
\eqref{w:off:B-B}, \eqref{w:off:C-B} and \eqref{w:off:B-C}  with
$q=p=\tilde p$, $\theta_{1}= D/\tilde p$, $\theta_{2}=0 $ and $c=0$.
We interpolate (by the real method since  we allow sublinear
operators) these estimates  with the ones coming from
$L^q(\mu)-L^r(\mu)$ off-diagonal estimates on balls. Thus $(q_\theta,
r_\theta)\in  {\mathcal{O}}(\T)$ where $ 1/q_{\theta}= \theta/\tilde
p + (1-\theta)/q$, $1/r_{\theta}=\theta/\tilde p+(1-\theta)/r$ and
$\theta \in (0,1)$. Choosing $\theta$ such that $q_{\theta}=p$ proves
that $(p, q(p))\in {\mathcal{O}}(\T)$ with $q(p)=r_{\theta}>p$.

In general, $\C(\T)$ has no further structure.  For example on
$\RR^n$ equipped with Lebesgue measure, let $T_{t}$  be the operator
of convolution   with $t^{-n/2}\phi(x/\sqrt t)$ with $t>0$,  $\phi$
positive, supported in the unit ball and $\phi \in L^s$ if and only
if $1\le s\le \rho$ for some $\rho\in (1,\infty)$. From Young's
inequality, it is easy  to determine $\C(\T)$ as the region in
$[0,1]^2$ below the diagonal $v=u$ (excluded) and above the line
$v=u-1/\rho'$ (included) and also to find that
$\widetilde\J(\T)=[1,\infty]$.  In particular, there is no interval
$I$ in $[1,\infty]$ such that  for all $p,q$ with $p< q$, $p,q \in I$
is equivalent to  $T_{t}\in \off{p}{q}$,  in such a case $\C(\T)$
(and ${\mathcal{O}}(\T)$) would be a triangle.

\subsection{Application to semigroups}\label{sec:semigroups}

Let  $\B$ be a Banach space of measurable functions stable under
truncations with indicator functions of measurable sets and
containing all simple functions. In this way, $\B\cap L^p(\mu)$ is
dense in $L^p(\mu)$ for all $1\le p<\infty$.

Let $\{T_{t}\}_{t>0}$ be a semigroup of bounded linear operators on
$\B$, that is, we assume for $t,s>0$ that
$$
T_t \in \call (\B);
\qquad\qquad
T_s\circ T_{t}=T_{s+t}.
$$
Here and in what follows $\call (X)$ denotes the set of
bounded linear operators on a Banach space  $X$.

\begin{prop}\label{prop:int-J} Set $\T=\{T_{t}\}_{t>0}$.  Assume there exist $\tilde p, \tilde q$ with
$1 \le \tilde p < \tilde q \le \infty$ such that $T_{t}\in
\off{\tilde p}{\tilde q}$.\footnote[2]{It is understood that the
functions to be considered are in $L^{\tilde p}(\mu)\cap \B$.} Then,
there exists a unique subset of $[1,\infty]$, which we denote by
$\J(\T)$, such that the following holds\textup{:}
\begin{equation}
\label{eq:J-off}
\forall\, p,q \in [1,\infty], \ p< q \quad \big(\
T_{t}\in \off{p}{q} \Longleftrightarrow p,q \in \J(\T) \ \big).
\end{equation}
This set is an interval,  contains $[\tilde p, \tilde q]$,
$\J(\T)\subset \widetilde\J(\T)$ and
$\Int\J(\T)=\Int\widetilde\J(\T)$.
\end{prop}

\begin{remark}\rm
This propagation property is reminiscent of the extrapolation for $L^p-L^q$ boundedness developed for semigroups in \cite{Cou}.
\end{remark}

With the notation of the previous section, \eqref{eq:J-off}
reformulates into
$$
\forall\, p,q \in [1,\infty], \ p< q \quad \big( \ (p,q) \in
{\mathcal{O}}(\T) \Longleftrightarrow p,q \in \J(\T)\ \big),
$$
which means that ${\mathcal{O}}(\T)$ is a triangle.

\begin{proof}

Note that if $E,F$ are two subsets such that \eqref{eq:J-off} holds
for $E$ and $F$ then, clearly, $E=F$ and so the uniqueness follows.
Let us now construct such a set.

Fix $\tilde p < \tilde r < \tilde q$. Let $\J_{-}(\T)$ be the set of
all $p\in [1,\tilde r]$ such that    $(p,\tilde r) \in
{\mathcal{O}}(\T)$. By one of the remarks after Definition
\ref{defi:off-d:weights}, this set is an interval  with upper bound
$\tilde r$ and it contains $[\tilde p, \tilde r]$. Similarly, the set
$\J_{+}(\T)$  of  all $p\in [\tilde r, \infty]$ such that  $(\tilde
r, p) \in {\mathcal{O}}(\T)$,   is an interval containing $[\tilde r,
\tilde q]$. Set $\J(\T)= \J_{-}(\T)\cup \J_{+}(\T)$. This  is clearly
an interval and it contains $[\tilde p, \tilde q]$.

Let us see that $p,q \in \J(\T)$ with $p<q$ imply $T_{t}\in
\off{p}{q}$. Indeed, if $p<q\le \tilde r$ or $\tilde r \le p < q $,
then $T_{t}\in \off {p}{q}$ using one of the remarks after Definition
\ref{defi:off-d:weights}, hence $(p,q) \in  {\mathcal{O}}(\T)$. If $p
\le \tilde r < q$, then $T_{t}\in \off {p}{\tilde r}$ and $T_{t}\in
\off {\tilde r}{q}$ hence by the semigroup property and Theorem
\ref{theor:uniform-compo}, part $(b)$,  $T_{2t}\in \off {p}{q}$. But
we may change $2t$ to $t$ and we have $(p,q) \in  {\mathcal{O}}(\T)$.

We prove  the converse: let $1\le p < q \le \infty$ with $T_{t}\in
\off{p}{q}$ and let us show that $p,q \in \J(\T)$.

\

\noindent \textbf{Case $q\le \tilde r$:} We have $T_{t}\in
\off{p}{q}$ and $T_{t}\in \off{\tilde r}{\tilde r}$.  Hence, by
interpolation, $T_{t}\in \off{p_\theta} {q_\theta}$ where
$1/p _{\theta}= \theta/ p + (1-\theta)/\tilde r$,
$1/q _{\theta}=\theta/q+(1-\theta)/\tilde r$ and $\theta \in (0,1)$.
If $p<\inf \J_{-}(\T)$ then we can choose $\theta$ such that
$p_{\theta}<\inf \J_{-}(\T)<q_{\theta}$. Since $\J_{-}(\T)$ is an
interval,  $q_{\theta}\in \J_{-}(\T)$, that is, $T_{t}\in
\off{q_{\theta}}{\tilde r}$. By Theorem \ref{theor:uniform-compo},
part $(b)$, and the semigroup property, $T_{2t}\in \off
{p_{\theta}}{\tilde r}$. Changing $2t$ to $t$ proves that $p_{\theta}
\in \J_{-}(\T)$, which is a contradiction. We have therefore shown
that $p\ge \inf \J_{-}(\T)$. If $p>\inf\J_{-}(\T)$, then $p\in
\J_{-}(\T)$ as $p< \tilde r$ and $\J_{-}(\T)$ is an interval. If
$p=\inf\J_{-}(\T)$, then $q \in \J_{-}(\T)$, hence $T_{t}\in
\off{q}{\tilde r}$. As $T_{t}\in \off{p}{q}$ by assumption, we have
again $T_{2t}\in \off{p}{\tilde r}$, hence $p\in \J_{-}(\T)$. We have
shown in this case that both $p$ and $q$ belong to $\J_{-}(\T)
\subset \J(\T)$.

\

\noindent \textbf{Case $p\ge \tilde r$:} This case is similar to the
previous one  by changing $\inf \J_{-}(\T)$ to $\sup \J_{+}(\T)$
(where the supremum is $\infty$ if $\J_{+}(\T)$ is unlimited) and
arguing on $q$ in place of $p$.

\

\noindent \textbf{Case $p<\tilde r<q$:}  By one of the remarks after
Definition \ref{defi:off-d:weights}, we have that $T_{t}\in
\off{p}{\tilde r}$ and $T_{t}\in \off{\tilde r}{q}$. Hence, by
definition, $p\in \J_{-}(\T)$ and $q \in \J_{+}(\T)$.

\

Let us finish the proof  by comparing the interiors of $\J(\T)$ and
$\widetilde \J(\T)$. By Theorem \ref{theor:uniform-compo}, $\J(\T)
\subset \widetilde \J(\T)$,  and the inclusion passes to interiors.
Since $\tilde p, \tilde q \in \J(\T)$, ${\mathcal{O}}(\T)\ne
\emptyset$. We showed in the previous section that  then for each
$p\in \Int\widetilde \J(\T)$, $T_{t} \in \off{p}{q}$ for some
$q=q(p)>p$. In particular, $p\in \J(\T)$ by \eqref{eq:J-off}. Thus,
$\Int\widetilde\J(\T) \subset \Int\J(\T)$.
\end{proof}

The following result shows that off-diagonal estimates on balls for a
semigroup   propagate to a sectorial analytic extension with
\textbf{optimal angle of sectors} provided there is one pair $(p_{0},
p_{0})$ for which one has off-diagonal estimates of balls  for the
analytic extension.

We consider $\{T_z\}_{z\in \Sigma_\vartheta}$  an analytic semigroup
of bounded linear operators on $\B$ with angle $\vartheta<\pi/2$,
that is, we assume for $z,z' \in \Sigma_{\vartheta}=\{\zeta\in
\mathbb{C}\setminus \{0\}\, ; \, |\arg \zeta| < \vartheta\}$,
$$
T_z \in \call (\B);
\qquad\qquad
T_z\circ T_{z'}=T_{z+z'};
\qquad\qquad
z\in \Sigma_{\vartheta} \longmapsto T_{z} \in \call (\B) \mbox{
is  analytic}.
$$
We say that $\{T_z\}_{z\in \Sigma_\vartheta}\in\off{p}{q}$ whenever
it satisfies the estimates in Definition \ref{defi:off-d:weights}
with $|z|$ in place of $t$. By density, this implies in particular
that the semigroup has an analytic extension from $\Sigma_{\vartheta}$ into $\call(L^r(\mu))$ for
$p\le r \le q$.

Recall that $\widetilde \J(\T)$ denotes the maximal interval of those  $p\in [1,\infty]$
for which $T_{t}$ is      bounded on $L^p(\mu)$ uniformly in $t>0$.

\begin{theor}\label{theor:off:analytic}
Let $1\le p\le p_0\le q\le \infty$ and $\vartheta_1$ with $0\le
\vartheta_1<\vartheta$. Assume that $\{T_t\}_{t>0}\in\off{p}{q}$ and
that $\{T_z\}_{z\in \Sigma_{\vartheta}}\in \off{p_0}{p_0}$. Then
for any $m\in \NN$, $\{z^m \frac {d^m T_z}{dz^m}\}_{z\in \Sigma_{\vartheta_1}}\in \off{p}{q}$.

\end{theor}

\begin{proof} Assume first  $m=0$.
Any $z\in \Sigma_{\vartheta_1}$ has a decomposition $z=s+w+t$ where
$w\in \Sigma_{\vartheta}$, $s,t >0$ and $|z|\sim |w| \sim s \sim t$,
the constants of comparability depending only on $\vartheta,
\vartheta_{1}$. Hence, we can write $T_{z}=T_{s}\circ T_{w}\circ
T_{t}$ and use $T_{t }\in \off{p}{p_{0}}$, $T_{w}\in
\off{p_{0}}{p_{0}} $ and $T _{s}\in \off{p_{0}}{q}$ together with
part $(b)$ in Theorem \ref{theor:uniform-compo}.

For $m>0$, we use a third angle $\vartheta_{2}$ with $
\vartheta_1<\vartheta_{2}<\vartheta$. We just showed that
$\{T_z\}_{z\in \Sigma_{\vartheta_{2}}}\in \off{p}{q}$. To conclude we
only have to use Cauchy formulae on circular on circular contours we
can compute $\frac {d^m T_z}{dz^m}$ for $z\in \Sigma_{\vartheta_{1}}$
from $T_{\zeta}$ with $\zeta\in \Sigma_{\vartheta_{2}}$.
\end{proof}

Note that the assumption on the analytic semigroup is
$\off{p_0}{p_0}$ for the same exponent $p_{0}$ at both places. In
applications, $p_{0}=2$ arises often (see the introduction) but with a weight this exponent is no longer natural.

\

So far, we were only concerned about the action of the semigroup
operator  $T_{t}$  on $L^p(\mu)$ and its off-diagonal estimates.
Recall that $\widetilde \J(\T)$ is the interval of exponents $p$
such that it has an extension to a bounded semigroup to
$L^p(\mu)$. To define an infinitesimal generator, it suffices that
(the extension to $L^p(\mu)$ of)  the semigroup is continuous at
$0$ for the strong topology in $\call(L^p(\mu))$. As usual, we
remove $p=\infty$ from the discussion. However, the off-diagonal
estimates play a crucial role.

\begin{prop}\label{prop:generator}
Assume that Proposition \ref{prop:int-J} applies and that there is
some $r\in \J(\T)$, $r\ne \infty$, such that   $T_{t}$ is strongly
continuous on $L^r(\mu)$. Then, $T_{t}$ is strongly continuous on
$L^p(\mu)$ for all $p\in \J(\T)$ with $p\ne \infty$.  In particular,
it has an infinitesimal generator on those $L^p(\mu)$.
\end{prop}

\begin{proof}
 If $p
\in \J(\T)$ with  $p<r$,
for $f$ any simple function supported in a ball $  B$  we deduce that
$$
\Big(\aver{2\,B} |T_{t}f - f|^p\,  d\mu\Big)^{\frac 1 p }
\le
\Big(\aver{2\,B} |T_{t} f - f|^r\,  d\mu\Big)^{\frac 1 r }  \longrightarrow 0
$$
as $t\to 0$.
Next, the off-diagonal estimates on balls imply that
$$
\Big(\int_{(2\, B)^c} |T_{t} f - f|^p\,  d\mu\Big)^{\frac 1 p }
=
\Big(\int_{(2\,B)^c} |T_{t} f |^p\,  d\mu\Big)^{\frac 1 p } \longrightarrow 0
$$
as $t\to 0$, using the support of $f$ and Lemma \ref{lemma:Bc} below.
Then a density argument shows the strong continuity in $L^p(\mu)$.

 If $p\in \J(\T)$, $r<p<\infty$,  then the above applies to the
dual semigroup and we can use the well-known fact that on a reflexive space, the dual semigroup of a strongly continuous bounded semigroup is also strongly continuous (see, \textit{e.g.} \cite[Chapter 1]{Da1}).
\end{proof}

Let us turn to weighted off-diagonal  estimates. Assume  that $\T$ is
a semigroup as in Proposition \ref{prop:int-J}. Let  $w\in A_\infty$.
As $(\X,d,w)$ is a space of homogeneous type, we can apply
Proposition \ref{prop:int-J} provided we have some off-diagonal
estimates to start with. In this case, we can define an interval
$\J_w(\T)$ characterized as the unique set $E$ in $[1,\infty]$ for
which  whenever
$ 1\le p<q\le \infty$ the property $T_{t}\in \offw{p}{q}$ is equivalent to  $p,q \in E$.  Also, $\widetilde{\J}_w(\T)$ is the interval
of those $p\in [1,\infty]$  for which $T_{t}$ is bounded uniformly in
$t$ on $L^p(w)$.

Given $1\le p_0< q_0\le\infty $  we define the set
$$
\W_w(p_0,q_0)
=
\big\{
p: p_0<p<q_0, w\in A_{\frac{p}{p_0}}\cap
RH_{\left(\frac{q_0}{p}\right)'}
\big\}.
$$

\begin{corol}\label{corol:off-w}
 Let $1\le p_0< q_0\le\infty$ be such that $(p_0,q_0)\subset \J(\T)$ and
assume that $\W_w(p_0,q_0)\neq\emptyset$. Then, $\W_w(p_0,q_0)\subset
\J_w(\T)\subset \widetilde \J_{w}(\T)$ and, consequently, $\Int \J_w(\T) =\Int \widetilde{\J}_w
(\T)$. If, furthermore, $\T$ is strongly continuous on $L^{r}(\mu)$ for some $r\in (p_{0},q_{0})$, then $\T$ has an infinitesimal generator in $L^p(w)$ for all $p\in \J_{w}(\T)$, $p\ne \infty$.

\end{corol}

\begin{proof}  The first statement  is a consequence of
Proposition \ref{prop:off-unw-w} and the second of  Proposition \ref{prop:int-J} in this context
together with the fact  shown in \cite{AM1} that if
$\W_w(p_0,q_0)\neq\emptyset$ then it is an open interval. Concerning the last statement, by Proposition \ref{prop:generator} it suffices to check that $\T$ is strongly continuous on $L^p(w)$ for one $p\in \J_{w}(\T)$.

Choose $p\in \W_w(p_0,q_0)$. Then there exists $p_1$ with
$p_{0}<p<p_{1}<q_{0}$ and $w\in RH_{(\frac {p_1}{p})'}$. Hence,
$$
\Big(\aver{B} g^{p}\,  dw\Big)^{\frac 1 p }  \lesssim \Big (\aver{B} g^{p_{1}} \, d\mu \Big)^{\frac 1 {p
_{1}}},
$$
for any ball $B$ and positive measurable function $g$.  If we apply
this to $g=|T_{t} f - f|$ for $f$ any simple function supported in a
ball $ B$ (the hypothesis contains the fact that $T_{t}$ is defined
on $L^{p_1}(\mu)$) and let $t \to 0$, we deduce that
$$
\Big(\int_{2\,B} |T_{t} f - f|^p\,  dw\Big)^{\frac 1 p }
\lesssim
\Big(\int_{2\, B} |T_{t} f - f|^{p_{1}}\,  d\mu\Big)^{\frac 1 p_{1} }  \longrightarrow
0,
$$
where we have used that $T_t$ is strongly continuous on
$L^{p_1}(\mu)$ by Proposition \ref{prop:generator}. Next, the
off-diagonal estimates on balls for $dw$ imply that
$$
\Big(\int_{(2\,B)^c} |T_{t} f - f|^p\,  dw\Big)^{\frac 1 p }
=
\Big(\int_{(2\, B)^c} |T_{t}f |^p\,  dw\Big)^{\frac 1 p } \longrightarrow 0, \quad t\to 0,
$$
using the support of $f$ and Lemma \ref{lemma:Bc} below. Then a
density argument shows the strong continuity in $L^p(w)$.
\end{proof}

\begin{remark}\rm  Note that to define $\J_{w}(\T)$, we only need the existence of some
pair $(p,q)$ with $p<q$ such that $T_{t}\in \offw{p}{q}$. Our statement here is a concrete realization of this assumption.
\end{remark}

\section{A case study}\label{sec:casestudy}

We work in the Euclidean space with the Lebesgue measure. Let
$A=A(x)$ be an $n\times n$ matrix of complex and $L^\infty$-valued
coefficients defined on $\re^n$. We assume that this matrix satisfies
the following ellipticity (or \lq\lq accretivity\rq\rq) condition:
there exist $0<\lambda\le\Lambda<\infty$ such that
$$
\lambda\,|\xi|^2
\le
\Re A(x)\,\xi\cdot\bar{\xi}
\quad\qquad\mbox{and}\qquad\quad
|A(x)\,\xi\cdot \bar{\zeta}|
\le
\Lambda\,|\xi|\,|\zeta|,
$$
for all $\xi,\zeta\in\co^n$ and almost every $x\in \re^n$. We have
used the notation
$\xi\cdot\zeta=\xi_1\,\zeta_1+\cdots+\xi_n\,\zeta_n$ and therefore
$\xi\cdot\bar{\zeta}$ is the usual inner product in $\co^n$. Note
that then
$A(x)\,\xi\cdot\bar{\zeta}=\sum_{j,k}a_{j,k}(x)\,\xi_k\,\bar{\zeta_j}$.
Associated with this matrix we define the second order divergence
form operator
$$
L f
=
-\div(A\,\nabla f),
$$
which is understood in the standard weak sense by means of a
sesquilinear form.

The operator $-L$ generates a $C^0$-semigroup
$\{e^{-t\,L}\}_{t>0}$ of contractions on $L^2$. We wish to
study weighted off-diagonal estimates for $\{e^{-t\,L}\}_{t>0}$
and $\{\sqrt t\, \nabla e^{-t\,L}\}_{t>0}$. Before we do so, we
recall what is known on unweighted off-diagonal estimates and give
some complements.

\begin{remark}\label{remark:full-off}
\rm Let us emphasize that  on $\RR^n$, full off-diagonal estimates
are equivalent to off-diagonal estimates on balls in the
unweighted situation by Proposition \ref{prop:full-off}. This
implies in particular that if $e^{-t\,L}\in \fullx{p}{q}$ for some
$1\le p< q\le \infty$ then (passing to off-diagonal on balls and
then going back to full off-diagonal estimates) it follows that
$e^{-t\,L}\in \fullx{p_1}{q_1}$ for all $p\le p_1\le q_1\le q$. We
will use this fact later.
\end{remark}

\subsection{The intervals $\J(L)$ and $\K(L)$}

 Define  $\widetilde \J(L)$  (we change slightly the previous  notation to
emphasize the dependence on $L$)
 as the
interval of those exponents $p\in [1,\infty]$ such that
$\{e^{-t\,L}\}_{t>0}$ is bounded in  $\call(L^p)$.

An almost complete study of $L^p-L^q$ full off-diagonal estimates with $p<q$ for the semigroup has been done in
\cite{Aus} (and the exponent $\theta$ of the definition must be $\frac
1 2 \big( \frac n p - \frac n q\big)$). According to Proposition
\ref{prop:int-J} and Proposition \ref{prop:full-off}, we have the
following result.

\begin{prop} There exists a unique subset of $[1,\infty]$,  denoted by $\J(L)$,
which is a non empty interval, such that
\begin{equation}
\label{eq:J} \forall\, p,q \in [1,\infty], \ p< q \quad \big(\
e^{-t\, L}\in \fullx{p}{q} \Longleftrightarrow p,q \in \J(L) \
\big).
\end{equation}
  Furthermore, $\J(L)\subset \widetilde \J(L)$ and
$\Int\J(L)=\Int\widetilde\J(L)$.
\end{prop}

See Section \ref{sec:propagation} for the meaning of ``interior.''
Write $p_{-}(L)$ and $p_{+}(L)$ as the lower and upper bounds in
$[1,\infty]$ of $\J(L)$. According to the results proved or cited in
\cite{Aus},
$$
\J(L)=\widetilde \J(L)= [1,\infty],
\qquad
\mbox{if } n=1,2,
$$
$$
p_{-}(L) < \frac{2n}{n+2} \text{\quad and\quad}
p_{+}(L)>\frac{2n}{n-2},
\qquad
\mbox{if }n\ge 3.
$$
Note that in dimensions $n\ge 3$, it is not clear what happens at the
endpoints for either  boundedness  or off-diagonal estimates: can one
have boundedness and no off-diagonal estimates? Is $\J(L)$ open in
$[1,\infty]$?

\

Let us turn to the gradient of the semigroup. Define   $\widetilde
\K(L)$ as the interval of those exponents $p\in [1,\infty]$ such
that
 $\{\sqrt t\,
\nabla e^{-t\,L}\}_{t>0}$ is bounded in $\call(L^p)$.  This set
has been studied in
\cite{Aus}. It is an interval in $[1,\infty]$. If $q_{-}(L)$ and $q_{+}(L)$
denote respectively its lower and upper bounds,
then it is shown that $q_{-}(L)=p_{-}(L)$ and $p_+(L)\ge
(q_+(L))^*$ where, given $q$, its Sobolev exponent $q^*$ is
defined as $q^*=n\,q/(n-q)$ if $q<n$ and $q^*=\infty$ otherwise.
Also, we always have $q_{+}(L)>2$ with $q_{+}(L)=\infty$ if $n=1$.

This was proved with the help of   full off-diagonal estimates.
Define $\K_{-}(L)$  as the set of all $p\in [1,2]$ such that
$\{\sqrt t \, \nabla e^{-t\,L}\}_{t>0}$ satisfies $L^{p}-L^{2}$
full off-diagonal estimates and $\K_{+}(L)$  be the set of all
$p\in [2,\infty]$ such that    $\{\sqrt t \, \nabla
e^{-t\,L}\}_{t>0}$ satisfies  $L^{2}-L^{p} $ full off-diagonal
estimates. Set $\K(L)= \K_{-}(L)\cup \K_{+}(L)$. This  is an
interval  by interpolation since $2\in \K(L)$ and it is shown in
\cite{Aus} that $\Int\K(L)=\Int\widetilde \K(L)$.  If $n=1$,
$\K(L)=[1,\infty]$ (see \cite{AMcT}).

We wish to give some further observations, not noticed in
\cite{Aus}, especially concerning the endpoints of $\K(L)$.

\begin{lemma}\label{lemma:fode}
Let $1\le p<2$. The following assertions are equivalent:
\begin{list}{$(\theenumi)$}{\usecounter{enumi}\leftmargin=.8cm
\labelwidth=0.7cm\itemsep=0.2cm\topsep=.2cm
\renewcommand{\theenumi}{\alph{enumi}}}
\item $e^{-t\,L} \in \fullx{p}{2}$.

\item $\sqrt t \, \nabla e^{-t\,L}\in  \fullx{p}{2}$.

\item $t\,  L\, e^{-t\,L}\in \fullx{p}{2}$.

\end{list}

\end{lemma}

\begin{proof}
To prove that $(a)$ implies $(b)$, we observe that $\sqrt t \, \nabla
e^{-t\,L}\in \fullx{2}{2}$ because $2\in \K(L)$. Hence by composing
with $(a)$ and using the semigroup property, we obtain $(b)$.

Similarly, $\sqrt t \, e^{-t\,L}\div A \in \fullx{2}{2}$ because of
duality and $2 \in \K(L^*)$, and the fact that multiplication by
$A(x)$ is bounded on $L^2$. Hence, from $(b)$, it follows that $\sqrt
t \, e^{-t\,L}\, \div  A \circ \sqrt t \, \nabla e^{-t\,L}\in
\fullx{p}{2}$. This operator is nothing but $-t\,L\, e^{-2\,t\, L}$
and this proves $(c)$.

Let us assume $(c)$. Pick $E$,$F$ two closed sets, $f\in L^p\cap L^2$
with support in $E$ and $L^p$-norm  1 and $g \in L^2$ with support
in $F$ and $L^2$-norm 1. Setting $h(t)=\langle e^{-t\, L}f,
g\rangle$, it suffices to prove $|h(t)|
\lesssim t^{-\theta}\,
\expt{-\frac{c\,d^2(E,F)}{t}}$ with
$\theta=\frac12\,(\frac{n}{p}-\frac{n}{2})$. Observe that our
assumption says that $th'(t)$ has such a bound.

 First,  $\lim_{t\to\infty} h(t)=0$: this is a consequence of the
bounded holomorphic functional calculus for $L$ on $L^2$ since $z\mapsto
e^{-tz}$ converges to 0 uniformly on compact subsets of $\Re z>0$.
Hence, we can write $ h(t)= - \int_{t}^\infty h'(s) \, ds. $
Plugging the bound for $sh'(s)$ into this integral yields $|h(t)|
\lesssim t^{-\theta}$. This bound suffices when $d^2(E,F) \le t$.

The second case is when $0<t< d^2(E,F)$. In particular $E$ and $F$
are disjoint. Then, one has $\lim_{s\to 0} h(s)=\langle f,
g\rangle=0$. As $h(t)=\int_{0}^t
 h'(s)\, ds,$ the bound for $sh'(s)$ easily yields $|h(t)| \lesssim
t^{-\theta}\, \expt{-\frac{c\,d^2(E,F)}{t}}$.
\end{proof}

\begin{lemma}\label{lemma:fode1} Assume $n\ge 2$. Let $1\le p < q$ with
$q^*<\infty$. If $\sqrt t \, \nabla e^{-t\,L}\in \fullx{p}{q}$,
then $e^{-t\, L} \in \fullx{p}{q^*}$.
\end{lemma}

\begin{proof}  By Proposition \ref{prop:full-off}, it suffices to show that $e^{-t\, L} \in \offx{p}{q^*}$.
To this end, we shall need the following form of  Sobolev's
inequality.

\begin{lemma}\label{lemma:sobo}
Let $1\le q<n$. If $g\in L^{1}(\RR^n)$ with $\nabla g \in
L^q(\RR^n)$ then $g\in L^{q^*}(\RR^n)$ and there is a constant
$C>0$ such that   for any ball $B$,
$$\Big(\int_{\RR^n \setminus B} |g|^{q^*}\, dx \Big)^{\frac 1 {q^*}} \le C
\Big(\int_{\RR^n \setminus B} |\nabla g|^{q} \, dx \Big)^{\frac 1 q}.$$
\end{lemma}

The first part of the lemma is non classical but easy:  let
$\varphi_{j}$ be a smooth mollifying sequence  and set
$g_{j}=\varphi_{j}* g$. Then $g_{j}\in L^1\cap L^\infty(\RR^n)$ and
$\nabla g_{j }=\nabla g * \varphi_{j}\in L^q(\RR^n)$, so that  in particular $g_{j}\in
W^{1,q}(\RR^n)$. Thus Sobolev's inequality on $\RR^n$ applies to each
$g_{j}$ and yields
$$\Big(\int_{\RR^n} |g_{j}|^{q^*}\, dx \Big)^{\frac 1 {q^*}} \le C
\Big(\int_{\RR^n} |\nabla g_{j}|^{q} \, dx \Big)^{\frac 1 q} \le C \Big(\int_{\RR^n} |\nabla g|^{q} \, dx \Big)^{\frac 1 q}\, \|\varphi\|_{1}.$$
Of course, $C$ is independent of $j$.  The conclusion that $g\in L^{q^*}(\RR^n)$ follows by
applying Fatou's lemma  to a subsequence.

We next show the desired Sobolev estimate. It suffices to obtain the
desired inequality for $B$ being the unit ball, the general case
follows by a change of variable with no change on the constant.
Besides, it is enough to assume that $g\in C^1_{0}(\RR^n)$ by density
in $W^{1,q}(\re^n)$. Then for any $x\notin B$, one has
$g(x)=-\int_{0}^\infty\pm
\partial_{j}g(x\pm te_{j})\, dt$ where $e_{j}$ is any  vector of the
canonical basis  and the choice of signs depends on the location of
$x$: positive signs  when $x_j\ge 0$ and negative signs when $x_{j
}<0$. With this choice of signs, note that  if $x\notin B$ we have
for all $t\ge 0$, $|x\pm t\,e_j|\ge |x|\ge 1$ and so $x\pm
t\,e_j\notin B$. Hence, for all $j=1, \ldots, n$, $|g(x)|
\le \int_{-\infty}^{+\infty} |\nabla g (x+t e_{j})|
\bigchi_{\RR^n\setminus B} (x+t e_{j})\, dt$. From there, one can follow the
standard argument first with $q=1$ and then with other values of
$q$ (see, e.g.
\cite{Bre}).

We come back to Lemma \ref{lemma:fode1},  beginning  with the proof
of \eqref{w:off:B-C}   with respect to $dx$. Let $B$ be a ball, $r$
its radius and  $f\in C^\infty_{0}(\RR^n)$ with support in $B$. Let
$j\ge 2$.   Observe that $g=e^{-t\, L}f$ satisfy the hypotheses of
Lemma \ref{lemma:sobo}. Indeed, the full off-diagonal estimates on
$L^2$ and the support of $f$ imply that $\int_{\RR^n} |g(x)|^2\,
e^{\,c\, |x-x_{B}|^2 /t} \, dx<\infty$ for some $c>0$ where $x_{B}$
is the center of $B$. Hence $g\in L^1(\RR^n)$ from Cauchy-Schwarz
inequality. Furthermore,  $\nabla g \in L^q(\RR^n)$ by our
assumption. Thus, by Lemma \ref{lemma:sobo} and since $\sqrt t \,
\nabla e^{-t\,L}\in \fullx{p}{q}$ we have
\begin{eqnarray*}
\lefteqn{\hskip -1cm
\Big(\aver{C_j(B)} |e^{-t\,L}
f|^{q^*}\,dx\Big)^{\frac1{q^*}}
\lesssim
(2^j\,
r)^{-\frac{n}{q^*}}\, \Big(\int_{\RR^n \setminus 2^{j}\, B}
|\nabla e^{-t\,L} f|^{q}\,dx\Big)^{\frac1{q}}}
\\
&\lesssim&
(2^j\, r)^{-\frac{n}{q^*}}\,
 \sum_{l\ge j}
\Big(\int_{C_l(B)} |\nabla e^{-t\,L}
f|^{q}\,dx\Big)^{\frac1{q}}
\\
&\lesssim&
(2^j\, r)^{-\frac{n}{q^*}}\,
 \sum_{l\ge j}  t^{-\frac 1 2 (\frac n p - \frac n q ) - \frac 1
2}\,\expt{-\frac{c\,4^l\,r^2}{t}}  \,
\Big(\int_{B} |f|^{p}\,dx\Big)^{\frac1{p}}
\\
&\lesssim&
 2^{-j\frac{n}{q^*}}\,
  \sum_{l\ge j} 2^{-l(\frac n p - \frac n {q^*})}\,   \left(\frac{2^l \,
r}{\sqrt{t}}\right)^{ \frac n p - \frac n {q^*}}
\expt{-\frac{c\,4^l\,r^2}{t}} \,
\Big(\aver{B} |f|^{p}\,dx\Big)^{\frac1{p}}
\\
&\lesssim&
 \dec{\frac{2^j \, r}{\sqrt{t}}}^{ \frac n p - \frac n {q^*}}
\, \expt{-\frac{c\,4^j\, r^2}{t}}\,
\Big(\aver{B} |f|^{p}\,dx\Big)^{\frac1{p}}.
\end{eqnarray*}
Hence we obtain  \eqref{w:off:B-C} for $e^{-t\, L} \in
\offx{p}{q^*}$ with $\theta_{1}=0$ and $\theta_{2}= \frac n p -
\frac n {q^*}$.

The proof of \eqref{w:off:B-B} for $e^{-t\, L} \in \offx{p}{q^*}$
is similar using  Sobolev's inequality on $\RR^n$ only since we do
not need a Gaussian term and we obtain the same values for
$\theta_{1}$, $\theta_{2}$.

It  remains to see \eqref{w:off:C-B} for $e^{-t\, L} \in
\offx{p}{q^*}$.  Let $B$ be a ball, $r$ its radius,   $j\ge 2$ and $f
\in C_{0}^\infty(\RR^n)$ with $\supp f\subset C_{j}(B)$. Since
$C_{j}(B)=2^{j+1}\, B\setminus 2^j\, B$, we can cover $C_{j}(B)$ by a
finite number of balls $B_{j,k}$ with radii $\frac 5 8  \, 2^j\, r$
with centers at distance $\frac 3 2 \, 2^j\, r$ from the center of
$B$ and the number of balls is a dimensional constant independent of
$j$ and $B$.  It is enough to assume that $f$ is also supported in
one $B_{j,k}$.  Then   observe that $B$ is contained in $\RR^n
\setminus  2 B_{j,k}$, hence the preceding argument  changing $B$ to
$B_{j,k}$ yields \eqref{w:off:C-B} with $\theta_{1}$ and $\theta_{2}$
as above . Details are left to the reader.

In this way we have shown that $e^{-t\, L} \in \offx{p}{q^*}$ and
by Remark \ref{remark:full-off} this completes the proof.
\end{proof}

\begin{prop}\label{prop:int-K(L)full}  Assume $n\ge 2$. We have  $\K(L) \subset
\J(L)$ and $\K(L)$ is characterized by \begin{equation}
\label{eq:Kfull} \forall\, p,q \in [1,\infty], \ p< q \ \big(\
\sqrt t \, \nabla e^{-t\,L}\in \fullx{p}{q} \ \Longleftrightarrow
p,q \in \K(L) \ \big).
\end{equation}

\end{prop}

In \cite{Aus}, it is only shown that $\Int \K(L) \subset \J(L)$
and the characterization is not considered.

\begin{proof}  From Lemma \ref{lemma:fode}, $\K_{-}(L)=\J(L)\cap [1,2]$. From
$p_+(L)\ge (q_+(L))^*> q_+(L)$, we have $\K_{+}(L) \subset \J(L)$.
It follows that $\K(L) \subset \J(L)$.

Let us see \eqref{eq:Kfull}. Assume that $p,q \in \K(L)$ with
$p<q$. If $p<q\le  2$ or $2 \le p < q $, then $\sqrt t \, \nabla
e^{-t\,L}\in \fullx{p}{q}$ as a consequence of $p, q \in
\K_{-}(L)$ or $p,q \in \K_{+}(L)$ using the equivalence between
full off-diagonal estimates and off-diagonal estimates on balls
(see  Remark \ref{remark:full-off}). If $p \le 2 < q$, then
$\sqrt t \, \nabla e^{-t\,L}\in \fullx{2}{q}$ and $e^{-t\,L}\in
\fullx{p}{2}$ by Lemma \ref{lemma:fode}. Hence, by composition and
the semigroup property, $\sqrt t \, \nabla e^{-t\,L}\in
\fullx{p}{q}$.

\

We turn to the converse. Let $1\le p < q \le \infty$ with $\sqrt t
\, \nabla e^{-t\,L}\in \fullx{p}{q}$ and let us show that $p,q \in
\K(L)$.

\

\paragraph{\bf Case $2\le p <q$:}
We have $\sqrt t \, \nabla e^{-t\,L}\in \fullx{p}{q}$  and $\sqrt t
\, \nabla e^{-t\,L}\in \fullx{2}{2}$.  Hence, by interpolation,
$\sqrt t \, \nabla e^{-t\,L}\in \fullx{p_{\theta}}{q_{\theta}}$ where
$1/p _{\theta}= (1-\theta)/ p + \theta/2$, $1/q
_{\theta}=(1-\theta)/q+\theta/2$ and $\theta \in (0,1)$.  If $p
\notin \K_{+}(L)$ then $q>\sup \K_{+}(L)$.  We can choose $\theta$
such that $p_{\theta}<\sup\K_{+}(L)<q_{\theta}$. Since
$\K_{+}(L)\subset \J(L)$, one has   $p_{\theta}\in \J(L)$, that is,
$e^{-t\,L}\in \fullx{2}{p_{\theta}}$. By composition and the
semigroup property,  $\sqrt t \, \nabla e^{-t\,L}\in
\fullx{2}{q_{\theta}}$, hence $q_{\theta} \in \K_{+}(L)$. This is a
contradiction. We have therefore shown that  $p\in \K_{+}(L)$. As we
have $\sqrt t \, \nabla e^{-t\,L}\in \fullx{p}{q}$ by assumption  and
$e^{-t\,L}\in \fullx{2}{p}$ since $p\in \J(L)$, by composition and
the semigroup property,
 $\sqrt t \, \nabla e^{-t\,L}\in \fullx{2}{q}$. Hence
 $q\in \K_{+}(L)$.

 \

\paragraph{\bf Case $p<2 \le q$:} Since $\sqrt t \, \nabla e^{-t\,L}\in
\fullx{p}{q}$,  using the equivalence between off-diagonal estimates
on balls and full off-diagonal estimates (see Remark
\ref{remark:full-off}),  we have that $\sqrt t \, \nabla e^{-t\,L}\in
\fullx{2}{q}$  and $\sqrt t \, \nabla e^{-t\,L}\in \fullx{p}{2}$.
Hence, $p\in \K_{-}(L)$ and $q\in \K_{+}(L)$.

\

\paragraph{\bf Case $p<q< 2$:}  As $n\ge 2$, we have $q^*<\infty$. Hence, Lemma
\ref{lemma:fode1} yields in particular $p\in \J(L)$. As $p<2$, we
have $p\in \K_{-}(L)$ by Lemma \ref{lemma:fode} and since $p<q<2$,
$q\in \K_{-}(L)$ as well.
\end{proof}

Let us finish this section with analyticity issues. For $L$ as above,
there exists $\vartheta\in[0,\pi/2)$ depending only on the
ellipticity constants such that for all $f\in\mathcal{D}(L)$
$$
\big|\arg \langle  Lf,f\rangle\big|
\le
\vartheta.
$$
We take the smallest $\vartheta$ such that this estimate holds. In
this case, one can obtain that $L$ is of type $\vartheta$ and its
semigroup $\{e^{-t\,L}\}_{t>0}$ has an analytic extension to a
complex semigroup $\{e^{-z\,L}\}_{z\in\Sigma_{\pi/2- \vartheta}}$
of contractions on $L^2$.

Applying Theorem \ref{theor:off:analytic} with $p_0=2$ and
Proposition \ref{prop:full-off}, one can obtain
 full off-diagonal estimates for the family $\{(zL)^me^{-z\,L}\}_{z\in
\Sigma_{\mu}}$ in the range $\J(L)$ and a similar type of
arguments yields the same thing  for the family $\{\sqrt z \,
\nabla (zL)^me^{-z\,L}\}_{z\in \Sigma_{\mu}}$ in the range
$\K(L)$, where $0<\mu<\pi/2 - \vartheta$ and $|z|$ replaces $t$ in
the estimates. We skip details.

We gather here  a particular case for later use in \cite{AM3}.
Recall that $\Int\J(L)=\big(p_{-}(L), p_{+}(L)\big)$ and
$\Int\K(L)
=\big(q_{-}(L), q_{+}(L)\big)$.

\begin{prop} \label{prop:sg}
Fix $m\in \NN$ and $0<\mu <\pi/2-\vartheta$.
\begin{list}{$(\theenumi)$}{\usecounter{enumi}\leftmargin=.8cm
\labelwidth=0.7cm\itemsep=0.3cm\topsep=.3cm
\renewcommand{\theenumi}{\alph{enumi}}}
\item If $p, q \in
\big(\,p_-(L),p_+(L)\big)$ with $p\le q$, then
$\{(zL)^me^{-z\,L}\}_{z\in \Sigma_{\mu}}$ satisfies $L^p-L^q$ full
off-diagonal estimates and is a  bounded set in $\call(L^p)$.

\item If $p,q \in
\big(\,q_-(L),q_+(L)\big)$ with $p\le q$, then $\{\sqrt z \,
\nabla (zL)^me^{-z\,L}\}_{z\in \Sigma_{\mu}}$ satisfies $L^p-L^q$
full off-diagonal estimates and is a bounded set in $\call(L^p)$.
\end{list}
\end{prop}

\subsection{The intervals $\J_{w}(L)$ and $\K_{w}(L)$}

As a consequence of Proposition \ref{prop:sg} and Proposition
\ref{prop:off-unw-w} we have the following result.

\begin{prop} \label{prop:sg-w}
Fix $m\in \NN$ and $0<\mu <\pi/2-\vartheta$.  Let $w \in
A_{\infty}$.
\begin{list}{$(\theenumi)$}{\usecounter{enumi}\leftmargin=.8cm
\labelwidth=0.7cm\itemsep=0.3cm\topsep=.3cm
\renewcommand{\theenumi}{\alph{enumi}}}
\item If $p, q \in \W_w\big(\,p_-(L),p_+(L)\big)$ with $p\le q$,
then $\{(zL)^me^{-z\,L}\}_{z\in \Sigma_{\mu}}$ satisfies
$L^p(w)-L^q(w)$ off-diagonal estimates on balls and is a bounded
set in $\call(L^p(w))$. \item If $p,q \in
\W_w\big(\,q_-(L),q_+(L)\big)$ with $p\le q$, then $\{\sqrt z \,
\nabla (zL)^me^{-z\,L}\}_{z\in \Sigma_{\mu}}$ satisfies
$L^p(w)-L^q(w)$ off-diagonal estimates on balls and is a bounded
set in $\call(L^p(w))$.
\end{list}
\end{prop}

This statement says that one has some  \textit{a priori} knowledge
of the intervals were we have weighted off-diagonal estimates on
balls. But, they could be larger than this. For a weight $w$, we
let $\widetilde \J_{w}(L)$ and $\widetilde \K_{w}(L)$ be the
intervals  of exponents $p\in[1,\infty]$  such that $e^{-t\, L}$
and $\sqrt t\, \nabla e^{-t\, L}$  respectively are   bounded on
$L^p(w)$ uniformly in $t>0$.

\begin{prop} \label{prop:sg-w:extension}
Fix $m\in \NN$ and $0<\mu <\pi/2-\vartheta$.  Let $w \in
A_{\infty}$.
\begin{list}{$(\theenumi)$}{\usecounter{enumi}\leftmargin=.8cm
\labelwidth=0.7cm\itemsep=0.3cm\topsep=.3cm
\renewcommand{\theenumi}{\alph{enumi}}}
\item  Assume $\W_{w}\big(p_{-}(L), p_{+}(L)\big)\ne \emptyset$.
There exists a unique subset of $[1,\infty]$,  denoted by
$\J_w(L)$,  which is an interval containing $\W_{w}\big(p_{-}(L),
p_{+}(L)\big)$, such that
\begin{equation}
\label{eq:Jw} \forall\, p,q \in [1,\infty], \ p< q \quad \big(\
e^{-t\, L}\in \offw{p}{q} \Longleftrightarrow p,q \in \J_{w}(L) \
\big).
\end{equation}
  Furthermore, $\J_{w}(L) \subset \widetilde\J_{w}(L)$ and  $\Int\J_{w}(L)=\Int\widetilde\J_{w}(L)$. Also
 if $p, q \in
\J_{w}(L)$ with $p\le q$, then $\{(zL)^me^{-z\,L}\}_{z\in
\Sigma_{\mu}}$ satisfies $L^p(w)-L^q(w)$ off-diagonal estimates on
balls and is a bounded set in $\call(L^p(w))$. \item Assume
$\W_{w}\big(q_{-}(L), q_{+}(L)\big)\ne \emptyset$. There exists a
subset of $[1,\infty]$, denoted by $\K_w(L)$, which is an interval
containing $\W_{w}\big(q_{-}(L), q_{+}(L)\big)$ with the following
properties: if $p,q \in \K_{w}(L)$ with $p\le q$ then  $\sqrt t\,
\nabla e^{-t\, L}\in \offw{p}{q} $ and, conversely, for $p,q \in
[1,\infty]$ with $p<q$ and $p\ne \inf \K_{w}(L)$, if $\sqrt t\,
\nabla e^{-t\, L}\in \offw{p}{q} $ then $p,q \in \K_{w}(L)$. In
particular, $\K_{w}(L) \setminus \{\inf\K_{w}(L)\}$ is the largest
open interval $I$ in $(1,\infty]$ characterized by
\begin{equation}
\label{eq:Kw} \forall\, p,q \in (1,\infty], \ p< q \quad \big(\
\sqrt t\, \nabla e^{-t\, L}\in \offw{p}{q} \Longleftrightarrow p,q
\in I \ \big).
\end{equation}
 Furthermore, $\K_{w}(L) \subset \widetilde\K_{w}(L)$ and $\Int\K_{w}(L)=\Int\widetilde\K_{w}(L)$. Also
 if $p, q \in
\K_{w}(L)$ with $p\le q$, then $\{ \sqrt z \, \nabla
(zL)^me^{-z\,L}\}_{z\in \Sigma_{\mu}}$ satisfies $L^p(w)-L^q(w)$
off-diagonal estimates on balls and is a bounded set in
$\call(L^p(w))$. \item  Let $n\ge 2$. Assume $\W_{w}\big(q_{-}(L),
q_{+}(L)\big)\ne \emptyset$. Then  $\K_{w}(L) \subset \J_{w}(L)$,
$\inf \J_{w}(L) = \inf \K_{w}(L)$ and $(\sup \K_{w}(L))^*_w \le
\sup \J_{w}(L) $. \item If $n=1$,  the intervals $\J_{w}(L)$ and
$\K_{w}(L)$ are the same and contain $(r_{w},\infty]$.
\end{list}
\end{prop}

We have set $q^*_w= \frac{q\, n \, r_{w}}{n\, r_{w} - q}$ when
$q<n\, r_{w}$ and $q^*_{w}=\infty$ otherwise. Recall that
$r_{w}=\inf\{r\ge 1\, ; \, w\in A_{r}\}$.

\begin{remark}\rm Let us assume that $L$ has real coefficients.
Then, the kernel of $e^{-t\, L}$ is
bounded above and below by Gaussians of the form $Ct^{-n/2}\,
\expt{-\frac{\alpha\,d^2(x,y)}{t}}$  with different constants in
each estimate.  Hence, for $p\ge 1$ and $w\in A_{\infty}$, we find
that $e^{-t\, L}$ is  bounded on $L^p(w)$ if and only if
$w\in A_{p}$. The sufficiency comes from the upper bound on the
kernel.  The necessity uses the positivity of the kernel and the
doubling condition on $w$ to derive  $w\in A_{p}$. Thus,
$\J_{w}(L)=\{p\in [1,\infty]\, : \, w\in A_{p}\}$. At the same time
$\W_{w}\big(p_{-}(L), p_{+}(L)\big)=\W_{w}(1,\infty)= (r_{w},\infty)$. If $w\in A_{1}$, then one has that
$\J_{w}(L)=[1,\infty]$. If $w\notin
A_{1}$,  $\J_{w}(L)=(r_{w},\infty]$. In all cases $\Int\J_{w}(L)=\W_{w}\big(p_{-}(L), p_{+}(L)\big)$.  The positivity of the semigroup  makes it in some
sense extremal among this class of  semigroups (for complex $L$).

Also $\K_{w}(L)=(r_{w}, k_{w}|$ where $k_{w}\ge
\frac{q_{+}(L)}{(s_{w})'} $ and $s_{w}=\sup \{ s \in (1,\infty]\, ;
\, w\in RH_{s}\}$ (whether $k_{w}$ is in $\K_{w}(L)$ is not known: we
suspect that $\K_{w}(L)$ is open in $[1,\infty]$).

This also shows that there is no upper bound of $\sup\J_{w}(L)$ in
terms of $\sup \K_{w}(L)$ as already observed for $w=1$.
\end{remark}

\begin{remark}\rm
We do not know examples where
  $\J_{w}(L)$ and $ \W_{w}\big(p_{-}(L), p_{+}(L)\big)$  have  different endpoints: such examples, if any,  must be complex.
\end{remark}

\begin{remark} \rm  It seems natural to expect that $\J_{w}(L)$ and $\K_{w}(L)$
are included in the set of $r\in [1,\infty]$ such that $w\in
A_{r}$. We are unable to show this.

Also, in part $(b)$, we lack of a general argument showing that if
$p=\inf\K_{w}(L)<q \le \inf \W_{w}\big( q_{-}(L), q_{+}(L)\big)=
q_{-}(L)r_{w}$ and $\sqrt t\, \nabla e^{-t\, L}\in \offw{p}{q} $ then
$p\in \K_{w}(L)$.

\end{remark}

\begin{proof}[Proof of Proposition \ref{prop:sg-w:extension}]
Part $(a)$ follows from Corollary \ref{corol:off-w}. For the statement corresponding to the family
$\{(z\,L)^m\,e^{-z\,L}\}$ we observe that, given $p,q\in\J_w(L)$ with $p\le q$,
there exists $p_0,q_0, r_0\in \J_w(L)$ so that $p_0\le p\le q\le
q_0$ and $p_0<r_0<q_0$ with $r_0\in \W_w\big(p_-(L),p_+(L)\big)$.
By using Proposition \ref{prop:sg-w} and Theorem
\ref{theor:off:analytic} it follows that
$\{(z\,L)^m\,e^{-z\,L}\}\in\offw{p_0}{q_0}\subset \offw{p}{q}$.

We next prove part $(d)$, part $(c)$ and part $(b)$ in this order. In fact, the construction of
$\K_{w}(L)$ is given during the proofs of part $(d)$ for dimension 1 and  part $(c)$ for higher dimensions.

\

\noindent {\em Proof of Proposition \ref{prop:sg-w:extension}, Part
$(d)$.} We recall that $p_{-}(L)=q_{-}(L)=1$ and
$p_{+}(L)=q_{+}(L)=\infty$, because the kernels of $e^{-t\, L} $ and
of $\sqrt t\,  \frac d{dx}e^{-t\, L}$ are pointwise dominated by
Gaussians [AMcT].   Hence, $\J_{w}(L) $ is the interval of those $p
\in [1,\infty]$ such that $e^{-t\, L} \in \offw{p}{\infty}$ and it
contains  $(r_{w},\infty]$. Define $\K_{w}(L) $ as the interval of
those $p \in [1,\infty]$ such that $\sqrt t\,  \frac d{dx}\,e^{-t\,
L} \in \offw{p}{\infty}$. We show that $\J_{w}(L)=\K_{w}(L)$.

 Let $p\in \J_{w}(L)$. As $\sqrt t\,  \frac d{dx}\,e^{-t\, L} \in
\offw{\infty}{\infty}$ and $e^{-t\, L} \in \offw{p}{\infty}$, we have
by composition and the semigroup property $\sqrt t\,  \frac
d{dx}\,e^{-t\, L} \in \offw{p}{\infty}$. Hence, $p\in \K_{w}(L)$.

Conversely, assume $\sqrt t\,  \frac d{dx}\,e^{-t\, L} \in
\offw{p}{\infty}$. Let $q\in \RR$ with $\frac q p  >r_{w}$ so that
$w\in A_{\frac q p}$. Let $B$ be a ball (an interval), $r$ its
radius and  $f \in C_{0}^\infty(B)$. Since $e^{-t\, L} f(x)$
vanishes at $\pm\infty$ by the compact support of $f$ and the
decay of the kernel of $e^{-t\, L} $, we have  for all $x\in \RR$,
$$
|e^{-t\, L} f(x)| \le  \int_{x}^\infty |(e^{-t\, L} f)'(y)| \, dy \le \int_{\RR} |(e^{-t\, L} f)'(y)| \, dy.
$$
Now,  with $C_{l}=C_{l}(B)$, we use $w\in A_{\frac{q}{p}}$ and our
assumption, which implies that $\sqrt t\,  \frac d{dx}e^{-t\, L}
\in \offw{p}{q}$,
\begin{eqnarray*}
\lefteqn{
\int_{\RR} |(e^{-t\, L} f)'(y)| \, dy
\lesssim \sum_{l\ge 1}
2^{l+1}\,r\, \aver{C_l}  |(e^{-t\, L} f)'(y)| \, dy}
\\
&\lesssim&
 \sum_{l\ge
1} 2^{l+1}\,r\,
\Big(\aver{C_l} |(e^{-t\, L} f)'(y)|^p \, dy\Big)^{\frac 1 p}
\\
&\lesssim&
 \sum_{l\ge
1} 2^{l}\,r\,
\Big(\aver{C_l} |(e^{-t\, L} f)'(y)|^q \, dw(y)\Big)^{\frac 1 q}
\\
&\lesssim&
r \, t^{-\frac 1 2}\,
\left(\dec{\frac{4\,r}{\sqrt{t}}}^{\theta_2}+ \sum_{l\ge 2}
2^{l\,(1+\theta_1)}\,\dec{\frac{2^l\,r}{\sqrt{t}}}^{\theta_2}\,
\expt{-\frac{c\,4^l\,r^2}{t}} \right)
\Big(\aver{B} |f|^{p}\,dw\Big)^{\frac1{p}}
\\
&\lesssim&
r \, t^{-\frac 1 2} \, \dec{\frac{r}{\sqrt{t}}}^{\max\{\theta_2,
1+\theta_1\}}\, \Big(1+\expt{-\frac{c\,r^2}{t}}\Big)\,
\Big(\aver{B} |f|^{p}\,dw\Big)^{\frac1{p}}
\\
&\lesssim&
\dec{\frac{ r}{\sqrt{t}}}^{\widetilde{\theta}_2}\,
\Big(\aver{B} |f|^{p}\,dw\Big)^{\frac1{p}}.
\end{eqnarray*}
In particular, this proves \eqref{w:off:B-B} for $e^{-t\, L} \in
\offw{p}{\infty}$.

Remark that if $x\in C_{j}$, then one has the more precise
estimate
$$
|e^{-t\, L} f(x)| \le \int_{\RR\setminus 2^{j}\, B} |(e^{-t\, L}
f)'(y)| \, dy.
$$
Indeed, it suffices to integrate $(e^{-t\, L} f)'$ from $x$ to
$+\infty$ if $x\ge 0$ and from  $-\infty$ to $x$  if $x\le 0$. In both
cases, the interval of integration is contained in $\RR\setminus
2^{j}\, B$. Hence, the same argument above yields
$$
|e^{-t\, L} f(x)|
\lesssim
\dec{\frac{ 2^j\, r}{\sqrt{t}}}^{\widetilde{\theta}_2} \,
\expt{-\frac{\alpha\,4^j\,r^2}{t}}\,\Big( \aver{B} |f|^p\,
dw\Big)^{ \frac 1 p},
$$
which proves \eqref{w:off:B-C} for $e^{-t\, L} \in
\offw{p}{\infty}$.

Similarly, assume $f$  supported in $C_{j}=C_{j}(B)$. We  want to
estimate $|e^{-t\, L} f(x)|$ for $x\in B$. Split $C_{j}$ into its
two connected components, $B_{j,1}, B_{j,2}$, which are intervals
of radius $2^{j-1}\, r$. Observe that $B$ is contained in
$\RR\setminus 2 \, B_{j,k}$ for $k=1,2$. Assume that $f$ is
supported in $B_{j,1}$ to fix ideas. Hence, for $x\in B$, one has
as before
$$
|e^{-t\, L} f(x)| \le \int_{\RR\setminus 2\, B_{j,1}} |(e^{-t\, L}
f)'(y)| \, dy.
$$
Arguing as above (with $2\, B_{j,1}$ in place of $2^j\,B$) we
obtain
$$
|e^{-t\, L} f(x)| \le \dec{\frac{2^j\,
r}{\sqrt{t}}}^{\widetilde{\theta}_2}\,
\expt{-\frac{\alpha\,4^j\,r^2}{t}}\,
\Big(\aver{B_{j,1}} |f|^{p}\,dw\Big)^{\frac1{p}}.
$$
One does the same thing when $f$ is supported in $B_{j,2}$. This
proves \eqref{w:off:C-B} for $e^{-t\, L} \in \offw{p}{\infty}$.
Hence, $p\in \J_{w}(L)$.

\

\paragraph{\em Proof of Proposition \ref{prop:sg-w:extension}, Part $(c)$}
We have $n\ge 2$, $\W_{w}\big(q_{-}(L), q_{+}(L)\big)\ne \emptyset$
and we know that this is an open interval. Pick $\tilde r \in
\W_{w}\big(q_{-}(L), q_{+}(L)\big)$ and set
\begin{align*}
\K_{-,w}(L)&=\{ p \in [1,\tilde r]\, ; \sqrt t \, \nabla\, e^{-t\, L}
\in \offw{p}{\tilde r}\},
\\[0.1cm]
\K_{+,w}(L)&=\{ p \in [\tilde r, \infty]\, ; \sqrt t\,  \nabla\,
e^{-t\, L} \in \offw{\tilde r}{p}\},
\\[0.1cm]
\K_{w}(L)&=\K_{-,w}(L) \cup \K_{+,w}(L).
\end{align*}
By construction, $\K_{w}(L)$ contains $ \W_{w}\big(q_{-}(L),
q_{+}(L)\big)$ and it is clearly an interval.

We need the following lemmas whose proofs are given below:

\begin{lemma}\label{lemma:K-wL} Let $1\le p <\tilde r$.
The following assertions are equivalent:
\begin{list}{$(\theenumi)$}{\usecounter{enumi}\leftmargin=.8cm
\labelwidth=0.7cm\itemsep=0.2cm\topsep=.2cm
\renewcommand{\theenumi}{\roman{enumi}}}
  \item $e^{-t\, L} \in \offw{p}{\tilde r}.$
  \item $\sqrt t \, \nabla \, e^{-t\, L} \in \offw{p}{\tilde r}.$
\end{list}
\end{lemma}

\begin{lemma}\label{lemma:K+wL} Assume $\tilde r < p \le \infty$ and
$\sqrt t \, \nabla \, e^{-t\, L} \in \offw{\tilde r}{p}$.  Then
for  $\tilde{r} \le q<p^*_{w}$, we have
 $e^{-t\, L} \in \offw{\tilde r}{q}.$
\end{lemma}

Note that Lemma \ref{lemma:K-wL} yields that $\inf \J_w(L)=\inf
\K_w(L)$ and $\K_{-,w}(L)\subset \J_{w}(L)$. On the other hand, Lemma \ref{lemma:K+wL} implies that
$(\sup \K_w(L))^*_w\le \sup \J_w(L)$ and so $\K_{+,w}(L)\subset
\J_w(L)$. This proves part $(c)$.

\begin{proof}[Proof of Lemma \ref{lemma:K-wL}]
As $\tilde r \in \W_{w}\big(q_{-}(L), q_{+}(L)\big)$, we have $\sqrt
t \, \nabla \, e^{-t\, L} \in \offw{\tilde r}{\tilde r}.$ Hence, by
composition and the semigroup property, we deduce that $(i)$ implies
$(ii)$.

For the converse, we cannot follow the route of Lemma
\ref{lemma:fode} so we use similar ideas as in Lemma
\ref{lemma:fode1}. We introduce some auxiliary exponents. Since
$\tilde r\in \W_w\big(q_-(L),q_+(L)\big)$, there exist $p_1$,
$q_1$ such that $q_-(L)<p_1<\tilde r <q_1<q_+(L)$ and $w\in
A_{\frac{\tilde r}{p_1}}\cap RH_{\left(\frac{q_1}{\tilde
r}\right)'}$. Note that $ q_+(L)> (q_{-}(L))^*$: indeed if $n=2$
then $q_-(L)=1$ and $q_+(L)>2$ whereas if $n\ge 3$, $q_-(L)<
\frac{2\,n}{n+2}$ and $q_+(L)>2$. Thus one can choose   $p_1<n$
and $q_{1}$ so  that $p_1^*=\frac{n\,p_1}{n-p_1}\le q_{1}$.

We begin  the proof  of $(i)$ with  \eqref{w:off:B-C}. Let $B$ be a
ball, $r$ its radius and  $f\in C^\infty_{0}(\RR^n)$ with support in
$B$. Let $j\ge 2$ and $C_{j}=C_{j}(B)$.   Observe that $g=e^{-t\,
L}f$ satisfy the hypotheses of  Lemma \ref{lemma:sobo} with
$q=p_{1}$. We use that $w\in RH_{\left(\frac{q_1}{\tilde r}\right)'}$
and $e^{-(t/2)\,L}\in \fullx{p_1^*}{q_1}$ ---because
$p_-(L)=q_-(L)<p_1<p_1^*\le q_1<q_+(L)\le p_+(L)$--- and Sobolev's
inequality on $\RR^n$ and on $\RR^n\setminus 2^{j-1}\, B $ (see Lemma
\ref{lemma:sobo}):
\begin{align*}
\lefteqn{\hskip-1cm
\Big(\aver{C_j} |e^{-t\,L}
f|^{\tilde r}\,dw\Big)^{\frac1{\tilde r}}
\lesssim
\Big(\aver{C_j} |e^{-t\,L} f|^{q_1}\,dx\Big)^{\frac1{q_1}}}
\\
&\lesssim
(2^j\, r)^{-\frac{n}{q_1}}\,
t^{-\frac12\,(\frac{n}{p_1^*}-\frac{n}{q_1})}\,
\bigg[ \expt{-\frac {c \, 4^j\, r^2}{t}}\,
\Big(\int_{2^{j-1}\, B} |e^{-(t/2)\,L} f|^{p_1^*}\,dx\Big)^{\frac1{p_1^*}}
\\
& \hskip6cm+
 \Big(\int_{\RR^n \setminus 2^{j-1}\, B} |e^{-(t/2)\,L}
f|^{p_1^*}\,dx\Big)^{\frac1{p_1^*}}\bigg]
\\
&\lesssim
(2^j\, r)^{-\frac{n}{q_1}}\,
t^{-\frac12\,(\frac{n}{p_1^*}-\frac{n}{q_1})}\,
\bigg[ \expt{-\frac {c \, 4^j\, r^2}{t}} \, \Big(\int_{\RR^n} |\nabla
e^{-(t/2)\,L} f|^{p_1}\,dx\Big)^{\frac1{p_1}}
\\
&\hskip6cm+
 \Big(\int_{\RR^n \setminus 2^{j-1}\, B} |\nabla e^{-(t/2)\,L}
f|^{p_1}\,dx\Big)^{\frac1{p_1}}\bigg].
\end{align*}
 From
$w\in A_{\frac{\tilde r}{p_1}}$ and  our assumption $(ii)$,  we have
\begin{eqnarray}
\lefteqn{\hskip-.5cm
\Big(\int_{\re^n} |\nabla e^{-(t/2)\,L}
f|^{p_1}\,dx\Big)^{\frac1{p_1}}
\lesssim \sum_{l\ge 1}
(2^{l+1}\,r)^{\frac n{p_1}}
\Big(\aver{C_l} |\nabla e^{-(t/2)\,L}
f|^{p_1}\,dx\Big)^{\frac1{p_1}}} \nonumber
\\
&\lesssim&
t^{-\frac12}\, \sum_{l\ge 1}(2^{l+1}\,r)^{\frac n{p_1}}\,
\Big(\aver{C_l} |\sqrt{t}\,\nabla e^{-(t/2)\,L}
f|^{\tilde r}\,dw\Big)^{\frac1{\tilde r}}
\nonumber
\\
&\lesssim&
r^{\frac n{p_1}} \, t^{-\frac 1 2}\,
\left(\dec{\frac{4\,r}{\sqrt{t}}}^{\theta_2}+ \sum_{l\ge 2}
2^{l\,(\frac{n}{p_1}+\theta_1)}\,\dec{\frac{2^l\,r}{\sqrt{t}}}^{\theta_2}\,
\expt{-\frac{c\,4^l\,r^2}{t}} \right)
\Big(\aver{B} |f|^{p}\,dw\Big)^{\frac1{p}}
\nonumber
\\
&\lesssim&
r^{\frac n{p_1}} \, t^{-\frac 1 2} \,
\dec{\frac{r}{\sqrt{t}}}^{\max\{\theta_2,
\frac{n}{p_1}+\theta_1\}}\, \Big(1+\expt{-\frac{c\,r^2}{t}}\Big)\,
\Big(\aver{B} |f|^{p}\,dw\Big)^{\frac1{p}}
\nonumber
\\
&\lesssim&
r^{\frac n{p_1}} \, t^{-\frac 1 2} \dec{\frac{
r}{\sqrt{t}}}^{\widetilde{\theta}_2}\,
\Big(\aver{B} |f|^{p}\,dw\Big)^{\frac1{p}}.
\label{aux-sob}
\end{eqnarray}

The integral  on $\RR^n\setminus 2^{j-1}B$ is analyzed similarly
when $j\ge3$ with a summation over $l\ge j-1$. If $j=2$, then the
integral  on $C_{1}=4B$ is replaced by one on $\widehat
C_{1}=4B\setminus 2B$ whose contribution is of the same order than
the one on $C_{2}$. Hence,
$$
 \Big(\int_{\RR^n \setminus 2^{j-1}\, B} |\nabla e^{-(t/2)\,L}
f|^{p_1}\,dx\Big)^{\frac1{p_1}}
 \lesssim (2^j\,r)^{\frac n{p_1}} \, t^{-\frac 1 2}\,
\dec{\frac{2^j \, r}{\sqrt{t}}}^{\widetilde{\theta}_2}\,
\expt{-\frac{c\, 4^j\,r^2}{t}}\,
\Big(\aver{B} |f|^{p}\,dw\Big)^{\frac1{p}}.
$$
All together after rearranging the terms yields
$$
\Big(\aver{C_j} |e^{-t\,L}
f|^{\tilde r}\,dw\Big)^{\frac1{\tilde r}}
\lesssim
2^{j\,\widetilde{\theta}_2}\, \dec{\frac{2^j \,
r}{\sqrt{t}}}^{\widetilde{\theta}_2 + \frac{n}{p_1}-\frac{n}{q_1} }\,
\expt{-\frac{c\, 4^j\,r^2}{t}}
\Big(\aver{B} |f|^{p}\,dw\Big)^{\frac1{p}},
$$
which is \eqref{w:off:B-C} for $e^{-t\, L} \in \offw{p}{\tilde
r}$.

The proof of \eqref{w:off:B-B} for $e^{-t\, L} \in \offw{p}{\tilde
r}$ is similar  using only Sobolev's inequality on $\RR^n$ since we
do not need a Gaussian term:
\begin{align*}
\Big(\aver{B} |e^{-t\,L}
f|^{\tilde r}\,dw\Big)^{\frac1{\tilde r}}
&\lesssim
\Big(\aver{B} |e^{-t\,L} f|^{q_1}\,dx\Big)^{\frac1{q_1}}
\lesssim
r^{-\frac{n}{q_1}}\, t^{-\frac12\,(\frac{n}{p_1^*}-\frac{n}{q_1})}\,
\Big(\int_{\re^n} |e^{-(t/2)\,L} f|^{p_1^*}\,dx\Big)^{\frac1{p_1^*}}
\\
&\lesssim
r^{-\frac{n}{q_1}}\, t^{-\frac12\,(\frac{n}{p_1^*}-\frac{n}{q_1})}\,
\Big(\int_{\RR^n} |\nabla
e^{-(t/2)\,L} f|^{p_1}\,dx\Big)^{\frac1{p_1}}.
\end{align*}
>From here we conclude the desired estimate as in \eqref{aux-sob}.

It  remains to prove \eqref{w:off:C-B} for $e^{-t\, L} \in
\offw{p}{\tilde r}$.
 Let $B$ be a ball, $r$ its radius,   $j\ge 2$ and $f \in C_{0}^\infty(\RR^n)$
with $\supp f\subset C_{j}=C_{j}(B)$. Since $C_{j}=2^{j+1}\,
B\setminus 2^j\, B$, we can cover $C_{j}$ by a finite number of
balls $B_{j,k}$ with radii $\frac 5 8  \, 2^j\, r$, with centers
at distance $\frac 3 2 \, 2^j\, r$ from the center of $B$ and the
number of balls is a dimensional constant independent of $j$ and
$B$. It is enough to assume that $f$ is also supported in one
$B_{j,k}$. Using  $w\in RH_{\left(\frac{q_1}{\tilde
r}\right)'}$ and $e^{-(t/2)\,L}\in \fullx{p_1^*}{q_1}$,
\begin{align*}
\lefteqn{\hskip-1cm
\Big(\aver{B} |e^{-t\,L}
f|^{\tilde r}\,dw\Big)^{\frac1{\tilde r}}
\lesssim
\Big(\aver{B} |e^{-t\,L} f|^{q_1}\,dx\Big)^{\frac1{q_1}}}
\\
&\lesssim
r^{-\frac{n}{q_1}}\,
t^{-\frac12\,(\frac{n}{p_1^*}-\frac{n}{q_1})}\,
\bigg[ \expt{-\frac {c \, 4^j\, r^2}{t}} \,
\Big(\int_{\RR^n\setminus 2^{j-2}\, B} |e^{-(t/2)\,L}
f|^{p_1^*}\,dx\Big)^{\frac1{p_1^*}}
\\
& \hskip6cm +
 \Big(\int_{2^{j-2}\, B} |e^{-(t/2)\,L}
f|^{p_1^*}\,dx\Big)^{\frac1{p_1^*}}\bigg].
\end{align*}
We use Sobolev's inequality in $\re^n$ and split $\RR^n$ according to
the sets  $C_{l}(2^{j+1}\, B)$, $l\ge 1$. Then  $w\in A_{\frac{\tilde
r}{p_1}}$ and   $(ii)$ yield as before
\begin{align*}
\lefteqn{\hskip-2cm
\Big(\int_{\RR^n\setminus 2^{j-2}\, B} |e^{-(t/2)\,L}
f|^{p_1^*}\,dx\Big)^{\frac1{p_1^*}}
\lesssim \Big(\int_{\RR^n} |\nabla e^{-(t/2)\,L}
f|^{p_1}\,dx\Big)^{\frac1{p_1}}}
\\
&\lesssim
(2^{j}\, r)^{\frac{n}{p_1}}\, t^{-\frac 1 2}
\,\dec{\frac{2^j\,r}{\sqrt{t}}}^{\widetilde{\theta}_2}\,
\Big(\aver{2^{j+1}\, B} |f|^{p}\,dw\Big)^{\frac1{p}}
\end{align*}
with $\widetilde \theta_{2}=\max\{\theta_{2}, \frac n {p_{1}} +
\theta_{1}\}$. Next, observe
 that $2^{j-2}\, B$ is contained in $\RR^n \setminus  2 B_{j,k}$, hence by Lemma
 \ref{lemma:sobo},
$$
\Big(\int_{2^{j-2}\, B} |e^{-(t/2)\,L} f|^{p_1^*}\,dx\Big)^{\frac1{p_1^*}}
\lesssim
\Big(\int_{\RR^n \setminus  2 B_{j,k}} |\nabla e^{-(t/2)\,L}
f|^{p_1}\,dx\Big)^{\frac1{p_1}}.
$$
Using  a splitting of $\RR^n \setminus  2 B_{j,k}$ with  the rings $
\, 2^{l+1}\, B_{j,k}
 \setminus   \, 2^l\, B_{j,k}$ for $l\ge 1$,  $w\in A_{\frac{\tilde r}{p_1}}$
and  $(ii)$  together with Lemma \ref{lemma:1-ann}  give us
$$
\Big(\int_{\RR^n \setminus  2 B_{j,k}} |\nabla e^{-(t/2)\,L}
f|^{p_1}\,dx\Big)^{\frac1{p_1}}
\lesssim
(2^{j}\, r)^{\frac{n}{p_1}}\, t^{-\frac 1 2}
\,\dec{\frac{2^j\,r}{\sqrt{t}}}^{\widetilde{\theta}_2}\,
\expt{-\frac {c \, 4^j\, r^2}{t}} \,
\Big(\aver{B_{j,k}} |f|^{p}\,dw\Big)^{\frac1{p}}.
$$
Gathering our estimates, we deduce that
$$
\Big(\aver{B} |e^{-t\,L}
f|^{\tilde r}\,dw\Big)^{\frac1{\tilde r}}
\lesssim
\dec{\frac{2^j\,r}{\sqrt{t}}}^{\widetilde{\theta}_2 +
\frac{n}{p_1}-\frac{n}{q_1} }\, \expt{-\frac {c \, 4^j\, r^2}{t}}
\,
\Big(\aver{C_{j}} |f|^{p}\,dw\Big)^{\frac1{p}}
$$
whenever $f$ is supported in $C_{j}\cap B_{j,k}$. This gives us
\eqref{w:off:C-B} for $e^{-t\, L} \in \offw{p}{\tilde r}$.
\end{proof}

\begin{proof}[Proof of Lemma \ref{lemma:K+wL}]
We know that $\tilde r \in \W_{w}\big(q_{-}(L), q_{+}(L)\big)$.
Hence, $w\in A_{\frac{\tilde r}{q_{-}(L)}} \subset A_{\tilde
r}\subset A_{p}$. Furthermore, for all $r>r_{w}$, all balls $B$
and Borel subsets $E$ of $B$,
$$
\frac{|E|}{|B|}\lesssim \bigg(
\frac{w(E)}{w(B)}\bigg)^{\frac{1}{r}}.
$$
Let $q<\infty$ with $\frac{1}{r}\ge \frac{n}{p}- \frac{n}{q}$.
Using
\cite[Corollary 3.2]{FPW}, we have an $L^p(w)-L^q(w)$ Poincar\'e inequality:
$$
\Big( \aver{B} \big|g-g_B\big|^q\, dw\Big)^{\frac{1}{q}} \lesssim
r_{B} \, \Big( \aver B |\nabla g|^p\, dw\Big)^{\frac{1}{p}},
$$
for all any $B$ and Lipschitz function $g$ where  $g_B$ stands for
the $w$-average of $g$ on $B$. Since convolution with a
$C_{0}^\infty$ function defines bounded map on $L^r(w)$ when $w\in
A_{r}$,   an approximation argument via mollifiers shows the validity
of this inequality if $g\in L^{\tilde r}(w)$ such that $\nabla g \in
L^p(w)$.

We begin with \eqref{w:off:B-B} for $e^{-t\, L} \in \offw{\tilde
r}{q}$. Since $\tilde r \in \W_{w}\big(q_{-}(L), q_{+}(L)\big) \subset
\W_{w}\big(p_{-}(L), p_{+}(L)\big)$,  $e^{-t\, L} \in \offw{\tilde
r}{\tilde r}.$ The matter is to improve integrability. Let $B$ be a ball, $r$ its radius and  $f\in
C^\infty_{0}(\RR^n)$ with support in $B$. Observe that the
Poincar\'e inequality above applies on $B$ to $g=e^{-t\, L}f$ since
we know that $g\in L^{\tilde r}(w)$,  $\nabla g \in L^p(w)$  from $e^{-t\, L} \in
\offw{\tilde r}{\tilde r}$,  our assumption $\sqrt t \, \nabla e^{-t\, L} \in
\offw{\tilde r}{p}  $ and Lemma \ref{lemma:Bc}. Hence
\begin{eqnarray*}
\Big(\aver B |e^{-t\, L} f|^q \, dw \Big)^{\frac{1}{q}}  &\lesssim& \aver B
|e^{-t\, L} f| \, dw  + r\, \Big(\aver B |\nabla e^{-t\, L} f|^p
\, dw \Big)^{\frac{1}{p}}.
\\
&\lesssim& \dec{ \frac r{\sqrt t}}^{\theta_{2}}  \, \Big(
\aver B |f|^{\tilde r}\, dw\Big)^{\frac{1}{\tilde r}}.
\end{eqnarray*}
 This proves
\eqref{w:off:B-B}.

 To prove \eqref{w:off:B-C}, we take $B$ and $f$ as before. Let  $j\ge 2$ and
cover  $C_{j}=C_{j}(B)$ by a finite number  of balls $B_{j,k}$
with radii $\frac58\, 2^j\, r$ and  centers at distance
$\frac32\,2^j\,r$ from the center of $B$. For each ball $B_{j,k}$,
we apply the same argument and obtain \eqref{w:off:B-C} using the
hypothesis with $C_{j}$ replaced by each $B_{j,k}$. It suffices to
add all the estimates to conclude.

 To prove \eqref{w:off:C-B}, we apply the same argument as for \eqref{w:off:B-B}
but with
 $f$ now supported in $C_{j}(B)$ for $j\ge 2$. Easy details are skipped.
\end{proof}

\paragraph{\em Proof of Proposition \ref{prop:sg-w:extension}, Part $(b)$} In parts $(c)$ and $(d)$, we defined a set
$\K_{w}(L)$ which is an interval in $[1,\infty]$ containing
$\W_{w}\big(q_{-}(L), q_{+}(L)\big)$.  The proof that  $p,q \in
\K_{w}(L)$ with $p\le q$ implies  $\sqrt t\, \nabla e^{-t\, L}\in
\offw{p}{q} $ is entirely similar to that of Proposition
\ref{prop:int-K(L)full} for $\K(L)$ replacing $2$ by $\tilde r$, full
off-diagonal estimates by off-diagonal estimates on balls and using
Lemma \ref{lemma:K-wL} in place of Lemma \ref{lemma:fode}.

Conversely assuming $p,q \in [1,\infty]$ with $p<q$ and $p\ne \inf
\K_{w}(L)$ and  $\sqrt t\, \nabla e^{-t\, L}\in \offw{p}{q} $ we
conclude that  $p,q \in \K_{w}(L)$ as in Proposition
\ref{prop:int-K(L)full} for $\K(L)$ except when $p<q\le \tilde r$.
For this situation, we argue as follows: As $p\ne \inf\K_{w}(L)$ we
have two cases. The first one is $p>\inf\K_{w}(L)$, which yields
$p,q$ in the interval $\K_{w}(L)$. The second one is $p<\inf
\K_{w}(L)$. Interpolating $\sqrt t\, \nabla e^{-t\, L}\in \offw{p}{q}
$ with $\sqrt t\, \nabla e^{-t\, L}\in \offw{\tilde p}{\tilde q} $
for any $\tilde p, \tilde q \in \K_{w}(L)$ with $\tilde p <\tilde q$,
one can find a pair $p_{\theta},q_{\theta}$ with $p_{\theta}
<\inf\K_{w}(L)$ and $q_{\theta }\in \W_{w}\big(q_{-}(L),
q_{+}(L)\big)$ such that $\sqrt t\, \nabla e^{-t\, L}\in
\offw{p_{\theta}}{q_{\theta}} $. Lemma \ref{lemma:K-wL} holds with
$q_\theta$ in place of $\tilde r$ as $q_\theta\in
\W_{w}\big(q_{-}(L), q_{+}(L)\big)$ and thus $p_\theta\in \J_w(L)$.
This leads to a contradiction since
$p_\theta<\inf\K_{w}(L)=\inf\J_{w}(L)$. Hence, this second case does
not happen.

\

The rest of the proof of $(b)$ is easy: That
$\Int\K_{w}(L)=\Int\widetilde \K_{w}(L)$ is a consequence of the
discussion in Section \ref{sec:propagation} and the previous
characterization (see also the last part of the proof of
Proposition \ref{prop:int-J}). The extension of the off-diagonal
estimates to the analytic family follows that done for the
semigroup. The $L^p(w)$ boundedness follows from Theorem
\ref{theor:uniform-compo}, part $(a)$. We skip further details.
\end{proof}

We conclude this discussion with a   word on infinitesimal generators.

\begin{corol}

Assume $\W_{w}\big(p_{-}(L), p_{+}(L)\big)\ne \emptyset$. For $p \in
\J_{w}(L)$ and $p\ne \infty$, the extension to $L^p(w)$
of  $\{e^{-t\, L}\}_{t>0}$  has an infinitesimal generator  which is an operator of type $\vartheta$ in $L^p(w)$.
\end{corol}

\begin{proof}  This is a consequence of Corollary  \ref{corol:off-w} and Proposition
\ref{prop:sg-w:extension}, noting that by construction $\{e^{-t\, L}\}_{t>0}$ is strongly continuous on $L^2$ and $2 \in (p_{-}(L), p_{+}(L))$. The fact that the infinitesimal generator is of type $\vartheta$ comes from the holomorphic extension of the semigroup on $\Sigma_{\pi/2-\vartheta}$ for $p\in \J_{w}(L)$.
\end{proof}

\begin{remark} \rm The inclusion $\K_{w}(L)\subset \J_{w}(L)$ implies that if $p\in \K_{w}(L)$ with $p\ne \infty$,  the domain of the $L^p(w)$-infinitesimal generator  is contained in the space
$\{f\in L^p(w)\, ; \, \nabla f \in L^p(w)\}$ (the gradient is defined in the distributional sense).
\end{remark}

\section{Proofs of the main results}\label{sec:proofs}

\subsection{Proof of Proposition \ref{prop:1infty}}
\label{sec:proof:prop:1infty}

Assume first that $K_{t}(x,y)$ is given with the desired properties.
Fix $t>0$. Let $B$ be a ball, $r$ its radius and $z$ its center.  Let
$f\in L^1(\mu)$ with support in $B$. Then for almost every $x\in B$,
$$
|T_{t}f(x)| \le \frac C {\mu(B(x,\sqrt t))}\, \int_{B}|f|\, d\mu \le
\frac{ C \, \mu(B)} {\mu(B(x,\sqrt t))}\,\Big(\aver{B}
|f|\,d\mu\Big).
$$
The doubling condition yields that $\mu(B)\approx \mu(B(x,r))$. If
$r\le \sqrt{t}$ then $\mu(B)\lesssim \mu(B(x,\sqrt{t})$. Otherwise
$r\ge \sqrt{t}$, the doubling condition implies
$$
\mu(B) \approx \mu (B(x,r))
\lesssim \Big(\frac{r}{\sqrt{t}}\Big)^D\,\mu (B(x,\sqrt{t})) ,
$$
and \eqref{w:off:B-B} holds with $\theta_{2}=D$, the doubling
exponent of $\mu$. Similarly, \eqref{w:off:C-B} and \eqref{w:off:B-C}
hold with $\theta_{1}=0$ and $\theta_{2}=D$. Hence $T_t
\in\off{1}{\infty}$.

Conversely, assume $T_t \in\off{1}{\infty}$. Fix $t>0$.  It follows in particular from \eqref{w:off:B-B} that for any ball $B$ and any $f,g \in L^1$ with support in $B$
$$
\int_{B}|g(x)|\, |T_{t}f(x)|\, d\mu(x) \le \frac C{\mu(B)}\,  \dec{\frac{r}{\sqrt{t}}}^{\theta_2} \|f\|_{1}\|g\|_{1}.$$
Hence, there exists $K_{t,B} \in  L^\infty(B\times B)$ such that
$$
\int_{B}g(x)\, T_{t}f(x)\, d\mu(x) = \int_{B}\int_{B} g(x) K_{t,B}(x,y)\, f(y) \, d\mu(y) \, d\mu(x).
$$
It is easy to show that $K_{t,B}(x,y)=K_{t,B'}(x,y) $ almost
everywhere on $B\times B \cap B'\times B'$ so that we may define
$K_{t}\in L^\infty_{\rm loc}(\X\times \X)$ which agrees almost
everywhere with $K_{t,B}$ on $B\times B$. Fix a Lebesgue point
$(x_{0}, y_{0})$ of $K_{t}$. Assume that $d(x_{0}, y_{0}) < \sqrt t$.
Fix $B=B(x_0,\sqrt{t})$ so that $x_{0},y_{0}\in B$. Then, apply the
formula above and let $f$, $g$ approximate Dirac masses at
$y_{0},x_{0}$ (more precisely, we use Lebesgue differentiation)  to
obtain
$$
|K_{t}(x_{0}, y_{0}) |\le  \frac C{\mu(B)} = \frac C {\mu(B(x_{0},
\sqrt t))}.
$$
If $d(x_{0},y_{0}) \ge \sqrt t$, then we choose $r= d(x_{0},
y_{0})/6$, and for $f \in L^1(\mu)$ with support in $B=B(y_{0}, r)$
and $g\in L^1(\mu)$ with support in $B'=B(x_{0},r)$, we have (embed
$B$ and $B'$ in a larger ball $B''$ on which the formula for
$K_{t,B''}$ is used)
$$
\int_{B'}g(x)\, T_{t}f(x)\, d\mu(x) = \int_{B'}\int_{B} g(x) K_{t}(x,y)\, f(y) \, d\mu(y) \, d\mu(x).
$$
Since $B'\subset C_{2}(B)$, we may apply  \eqref{w:off:B-C} with $j=2$ and by letting $f$ and $g$ approximate Dirac masses as before, we obtain
$$
|K_{t}(x_{0}, y_{0}) |\le  \frac C{\mu(B)} \, \Big(\frac{ d(x_{0},
y_{0})}{\sqrt t}\Big)^{\theta_{2}} \,   \expt{-\frac{c\, d^2(x_{0},
y_{0})}{t}}.
$$
But $d(x_{0},y_{0}) \ge \sqrt t$ implies that we can absorb the
$\theta_{2}$ power by the Gaussian factor and also that
$\mu(B(x_{0},\sqrt t)) \le  \mu(B(x_{0}, d(x_{0},y_{0}) ) )\lesssim
\mu(B(x_{0},r))=\mu(B)$ as $\mu$ is doubling.

\subsection{Proof of Theorem \ref{theor:uniform-compo}: Part $(a)$}
\label{sec:proof:theor:uniform-compo-a}

We need the following basic facts about spaces of homogeneous type.
Indeed, the following property was used originally to define those
spaces, see
\cite{CW}.

\begin{lemma}\label{lemma:geom:X}
There exists $N\in\NN$ depending on $C_0$ in \eqref{doubling}, such
that, for every $j\ge 1$, any ball $B$ contains at most $N^j$ points
$\{x_k\}_k$ such that $d(x_{k_1},x_{k_2})>r_B/2^{j}$.
\end{lemma}

We also recall the following well-known covering lemma whose proof is left to
the reader (note that the covering family has to be countable since
in any fixed ball the number of $r/2$-separated points is finite by
the previous result).

\begin{lemma}\label{lemma:covering-X}
Given $r>0$ there exists a sequence $\{x_k\}_{k}\subset\X$ so that
$d(x_{k_1},x_{k_2})>r/2$ for all $x_{k_1}\neq x_{k_2}$ and
$\X=\bigcup_k B(x_k,r)$.
\end{lemma}

We can now establish $(a)$ in Theorem \ref{theor:uniform-compo}. We
use Lemma \ref{lemma:covering-X} with $r=\sqrt{t}$ and write
$B_k=B(x_k,\sqrt{t})$. Then, if $1\le p<\infty$,
\begin{eqnarray*}
\|T_t f\|_{L^p(\mu)}^p
&\le&
\sum_{k} \int_{B_k}|T_t f|^p\, d\mu
\le
\sum_{k} \bigg( \sum_{j=1}^\infty
\Big(\int_{B_k}|T_t (\bigchi_{C_j(B_k)}\, f)|^p\, d\mu\Big)^{\frac 1 p }
\bigg)^{p}
\\
&\lesssim&
\sum_{k} \bigg( \sum_{j=1}^\infty
2^{j\,\theta_1}\,\dec{2^j}^{\theta_2}\, e^{-c\,4^j}\,
\left(\frac{\mu(B_k)}{\mu(2^{j+1}\,B_k)}\right)^{\frac 1 p }
\Big(\int_{C_j(B_k)}|f|^p\, d\mu\Big)^{\frac 1 p }
\bigg)^{p}
\\
&\lesssim&
\sum_{k} \bigg( \sum_{j=1}^\infty e^{-c\,4^j}\,
\int_{C_j(B_k)}|f|^p\, d\mu\bigg) \,
\bigg( \sum_{j=1}^\infty 2^{j\,(\theta_1+\theta_2)\, p'}\,e^{-c\,4^j}
\bigg)^{\frac p{p'}}\\
&\lesssim&
\sum_{j=1}^\infty e^{-c\,4^j}\, \int_{\X}|f|^p\, \sum_{k}
\bigchi_{C_j(B_k)}\,d\mu
\le
\sum_{j=1}^\infty e^{-c\,4^j}\, N^{j+3}\, \int_{\X}|f|^p\,d\mu
\\
&\lesssim& \int_{\X}|f|^p\,d\mu,
\end{eqnarray*}
where we have used that for any $j\ge 1$ we have $\sum_{k}
\bigchi_{C_j(B_k)}(x)\le N^{j+3}$. Indeed,  for a fixed $x\in\X$,
there exists $k_0$ such that $x\in B_{k_0}$. Then, by Lemma
\ref{lemma:geom:X},
$$
\sum_{k}  \bigchi_{C_j(B_k)}(x)
\le
\# \{k: x\in 2^{j+1}\,B_k\}
\le
\# \big\{k: x_k\in B(x_{k_0},2^{j+2}\,r)\big\}
\le
N^{j+3},
$$
since $d(x_k,x_j)>r/2=(2^{j+2}\,r)/2^{j+3}$. The modification for
$p=\infty$ is left to the reader.

\subsection{Proof of Theorem \ref{theor:uniform-compo}: Part $(b)$}
\label{sec:proof:theor:uniform-compo-b}

We next show that the definition of off-diagonal estimates on balls
is stable under composition.

To prove this  we need the following auxiliary results whose proofs
are postponed until the end of this subsection.

\begin{lemma}\label{lemma:sum}
Let $s>0$, $\alpha\ge 0$ and $\beta>0$ with $\alpha\neq \beta$. Then,
if $0<c'<c$,
$$
\sum_{k=0}^\infty
2^{k\,\alpha}\,\dec{2^k\,s}^{\beta}\,e^{-c\,4^k\,s^2}
\lesssim
\dec{s}^{\max\{\alpha,\beta\}}\,e^{-c'\,s^2}.
$$
\end{lemma}

\begin{remark}\label{remark:powers-diff}\rm
We have assumed $\alpha\neq \beta$ in order to get explicit exponents. If

 $\alpha=\beta$ the same estimate remains true with the
power of $\dec{s}$ being $\alpha+\varepsilon$, for any
$\varepsilon>0$, in place of $\alpha$. For cleanness and shortness,
we will use this lemma several times assuming that  the powers are
different, if this is not the case the final power has to be slightly
enlarged for the estimate to be correct.
\end{remark}

\begin{lemma}\label{lemma:1-ann}
If $T_t\in \off{p}{q}$ with exponents $\theta_1$ and $\theta_2$, then
for any ball $B$ with radius $r$ we have
$$
\Big(\aver{B}|T_t( \bigchi_{\widehat{C}_1(B)}\, f) |^{q}\,d\mu\Big)^{\frac 1 q}
\lesssim
\dec{\frac{2\,r}{\sqrt{t}}}^{\theta_2}\, \expt{-\frac{c\,4\,r^2}{t}}
\,
\Big(\aver{\widehat{C}_1(B)}|f|^{p}\,d\mu\Big)^{\frac 1 p }
$$
and
$$
\Big(\aver{\widehat{C}_1(B)}|T_t( \bigchi_B\, f) |^{q}\,d\mu\Big)^{\frac 1 q}
\lesssim
\dec{\frac{2\,r}{\sqrt{t}}}^{\theta_2}\, \expt{-\frac{c\,4\,r^2}{t}}
\,\Big(\aver{B}|f|^{p}\,d\mu\Big)^{\frac 1 p },
$$
that is, $T_t$ satisfies the last two estimates in Definition
\ref{defi:off-d:weights} with $j=1$ and $\widehat{C}_1(B)$ in place
of $C_j(B)$.
\end{lemma}

\begin{lemma}\label{lemma:Bc}
If $T_t\in \off{p}{q}$  with parameters $\theta_1, \theta_2, c$  then
for $0<c'<c$, for any ball $B$ with radius $r$ and for every $j\ge 1$
we have
$$
\Big(\aver{B}|T_t( \bigchi_{(2^j\,B)^c}\, f) |^{q}\,d\mu\Big)^{\frac 1 q}
\lesssim
2^{j\,\theta_1}\,
\dec{\frac{2^j\,r}{\sqrt{t}}}^{\max\{\theta_1,\theta_2\}}\,
\expt{-\frac{c'\,4^j\,r^2}{t}} \,
\Big(\aver{(2^j\,B)^c}|f|^{p}\,d\mu\Big)^{\frac 1 p }
$$
and
$$
\Big(\aver{(2^j\,B)^c}|T_t( \bigchi_B\, f)|^{q}\,d\mu\Big)^{\frac 1 q}
\lesssim
2^{j\,\theta_1}\,
\dec{\frac{2^j\,r}{\sqrt{t}}}^{\max\{\theta_1+D/q,\theta_2\}}\,
\expt{-\frac{c'\,4^j\,r^2}{t}}
\,\Big(\aver{B}|f|^{p}\,d\mu\Big)^{\frac 1 p }.
$$
\end{lemma}

\

Once we have stated these auxiliary results we can proceed to
establish Theorem \ref{theor:uniform-compo}, Part $(b)$.

\begin{proof}[Proof of Theorem \ref{theor:uniform-compo}: Part $(b)$]
We start with \eqref{w:off:B-B} and assume that $\supp f\subset B$.
Write $\lambda=r_B/\sqrt{t}$. Note that we have
\begin{eqnarray*}
\Big(\aver{B} |T_t(S_t f) |^{r}\,d\mu\Big)^{\frac1r}
&\le&
\Big(\aver{B} |T_t(\bigchi_{2\,B}\,S_t f) |^{r}\,d\mu\Big)^{\frac1r}
+
\Big(\aver{B} |T_t(\bigchi_{(2\,B)^c}\,S_t f) |^{r}\,d\mu\Big)^{\frac1r}
\\
&=&
I+II.
\end{eqnarray*}
Then, since $\supp f\subset B\subset 2\,B$ we have
\begin{eqnarray*}
I
&\lesssim&
\Big(\aver{2\,B} |T_t (\bigchi_{2\,B}\,S_t f) |^{r}\,d\mu\Big)^{\frac1r}
\lesssim
\dec{2\,\lambda}^{\theta_2}\,
\Big(\aver{2\,B} |S_t f|^{q}\,d\mu\Big)^{\frac1q}
\\
&\lesssim&
\dec{2\,\lambda}^{\theta_2+\gamma_2}\,
\Big(\aver{2\,B} |f|^{q}\,d\mu\Big)^{\frac1q}
\lesssim
\dec{\lambda}^{\theta_2+\gamma_2}\,
\Big(\aver{B} |f|^{p}\,d\mu\Big)^{\frac1p}.
\end{eqnarray*}
On the other hand, applying Lemma \ref{lemma:Bc} twice we have
\begin{eqnarray*}
II
&\lesssim&
2^{\theta_1}\,\dec{2\,\lambda}^{\max\{\theta_1,\theta_2\}}\,
e^{-c'\,4\,\lambda^2}\,
\Big(\aver{(2\,B)^c} |S_t f|^{q}\,d\mu\Big)^{\frac1q}
\\
&\lesssim&
2^{\theta_1+\gamma_1}\,
\dec{2\,\lambda}^{\max\{\theta_1,\theta_2\}+\max\{\gamma_1+D/q,\gamma_2\}}\,
e^{-c'\,4\,\lambda^2}\,
\Big(\aver{2\,B} |f|^{p}\,d\mu\Big)^{\frac1p}
\\
&\lesssim&
\dec{\lambda}^{\max\{\theta_1,\theta_2\}+\max\{\gamma_1+D/q,\gamma_2\}}\,\,
\Big(\aver{B} |f|^{p}\,d\mu\Big)^{\frac1p}.
\end{eqnarray*}
Collecting the bounds for $I$ and $II$ we obtain the desired
estimate.

Next, we consider \eqref{w:off:C-B}. Let $f$ be supported on $C_j(B)$
with $j\ge 2$. Let us set $\lambda=2^j\,r_B/\sqrt{t}$. We first split
the integral as follows
\begin{align*}
\Big(\aver{B} |T_t(S_t f) |^{r}\,d\mu\Big)^{\frac1r}
&\le
\Big(\aver{B} |T_t(\bigchi_{2^{j-1}\,B}\,S_t f) |^{r}\,d\mu\Big)^{\frac1r}
+
\Big(\aver{B} |T_t(\bigchi_{(2^{j-1}\,B)^c}\,S_t f) |^{r}\,d\mu\Big)^{\frac1r}
\\
&=
I+II.
\end{align*}
For $I$ we write $\widetilde{B}=2^{j-1}\,B$ which has radius
$r_{\widetilde{B}}=2^{j-1}\,r_B$. Thus,
$C_j(B)=4\,\widetilde{B}\setminus
2\,\widetilde{B}=\widehat{C}_1(\widetilde{B})$ and so by Lemma
\ref{lemma:1-ann} we have
\begin{align*}
I
&\le
\left(\frac{\mu(2^{j-1}\,B)}{\mu(B)}\right)^{\frac1r}
\Big(\aver{\widetilde{B}}
|T_t(\bigchi_{\widetilde{B}}\,S_t f) |^{r}\,d\mu\Big)^{\frac1r}
\\
&\lesssim
2^{j\,D/r}\, \dec{\frac{2\,r_{\widetilde{B}}}{\sqrt{t}}}^{\theta_2}\,
\,\Big(\aver{\widetilde{B}}|S_t(
\bigchi_{\widehat{C}_1(\widetilde{B})}\, f)|^{q}\,d\mu\Big)^{\frac1q}
\\
&\lesssim
2^{j\,D/r}\,\dec{\lambda}^{\theta_2}\,
\dec{\frac{r_{\widetilde{B}}}{\sqrt{t}}}^{\gamma_2}
\,\expt{-\frac{c\,4\,r_{\widetilde{B}}^2}{t}}
\Big(\aver{\widehat{C}_1(\widetilde{B})} |f|^{p}\,d\mu\Big)^{\frac1p}
\\
&\lesssim
2^{j\,D/r}\,\dec{\lambda}^{\theta_2+\gamma_2} \,e^{-c\,\lambda^2}\,
\Big(\aver{C_j(B)} |f|^{p}\,d\mu\Big)^{\frac1p}.
\end{align*}
On the other hand, by Lemma \ref{lemma:Bc}
\begin{eqnarray*}
II
&\lesssim&
2^{(j-1)\,\theta_1}\,
\dec{\frac{2^{j-1}\,r_B}{\sqrt{t}}}^{\max\{\theta_1,\theta_2\}}
\,\expt{-\frac{c'\,4^{j-1}\,r_B^2}{t}} \,
\Big(\aver{(2^{j-1}\,B)^c} |S_t f|^{q}\,d\mu\Big)^{\frac1q}
\\
&\lesssim&
2^{j\,\theta_1}\,\dec{\lambda}^{\max\{\theta_1,\theta_2\}}
\,e^{-c'\,\lambda^2}\,
\Big(\aver{(2^{j-1}\,B)^c} |S_t f|^{q}\,d\mu\Big)^{\frac1q}.
\end{eqnarray*}
Besides,
$$
\Big(\aver{(2^{j-1}\,B)^c} |S_t f|^{q}\,d\mu\Big)^{\frac1q}
\lesssim
\Big(\aver{2^{j+2}\,B}|S_t f|^{q}\,d\mu
\Big)^{\frac1q}
+
\Big(\aver{(2^{j+2}\,B)^c}|S_t f|^{q}\,d\mu
\Big)^{\frac1q}
=
II_1+II_2.
$$
For $II_1$ we observe that
\begin{eqnarray*}
II_1
&=&
\Big(\aver{2^{j+2}\,B}|S_t (\bigchi_{2^{j+2}\,B}\, f)|^{q}\,d\mu
\Big)^{\frac1q}
\lesssim
\dec{\frac{r_{2^{j+2}\,B}}{\sqrt{t}}}^{\gamma_2}\,
\,\Big(\aver{2^{j+2}\,B}|f|^{p}\,d\mu\Big)^{\frac1p}
\\
&\lesssim&
\dec{\lambda}^{\gamma_2}\,
\Big(\aver{C_j(B)}|f|^{p}\,d\mu\Big)^{\frac1p}.
\end{eqnarray*}
On the other hand for $II_2$ we set $\widetilde{B}=2^{j+1}\,B$ and so
its radius is $r_{\widetilde{B}}=2^{j+1}\,r_B$. Thus,
$(2^{j+2}\,B)^c=(2\,\widetilde{B})^c$ and $\,r_{\widetilde{B}}/{\sqrt t}= 2\lambda$.  By Lemma \ref{lemma:Bc}
we have
\begin{eqnarray*}
II_2
&=&
\Big(\aver{(2\,\widetilde{B})^c}
|S_t (f\,\bigchi_{\widetilde{B}}) |^{q}\,d\mu\Big)^{\frac1q}
\lesssim
2^{\gamma_1} \dec{4\, \lambda}^{\max\{\gamma_1+D/q,\gamma_2\}}\, e^{-
c'\,\lambda^2}\,
\Big(\aver{\widetilde{B}}|f|^{p}\,d\mu\Big)^{\frac1p}
\\
&\lesssim&
\dec{\lambda}^{\max\{\gamma_1+D/q,\gamma_2\}}\, e^{-c'\,\lambda^2}\,
\,\Big(\aver{C_j(B)}|f|^{p}\,d\mu\Big)^{\frac1p}.
\end{eqnarray*}
Collecting the bounds for $I$, $II_1$ and $II_2$ we
obtain the desired estimate.

Finally, we show \eqref{w:off:B-C}. We take $\supp f\subset B$ and
$j\ge 2$. Let us set $\lambda=2^j\,r_B/\sqrt{t}$. We proceed as
follows
\begin{eqnarray*}
\Big(\aver{C_j(B)} |T_t(S_t f) |^{r}\,d\mu\Big)^{\frac1r}
&\le&
\Big(\aver{C_j(B)} |T_t(\bigchi_{2^{j-1}\,B}\,S_t f) |^{r}\,d\mu\Big)^{\frac1r}
\\
&&
\hskip1.5cm +
\Big(\aver{C_j(B)} |T_t(\bigchi_{(2^{j-1}\,B)^c}\,S_t f) |^{r}\,d\mu\Big)^{\frac1r}
=
I+II.
\end{eqnarray*}
For $I$ we write $\widetilde{B}=2^{j-1}\,B$ which has radius
$r_{\widetilde{B}}=2^{j-1}\,r$. Thus,
$C_j(B)=4\,\widetilde{B}\setminus
2\,\widetilde{B}=\widehat{C}_1(\widetilde{B})$ and $r_{\widetilde{B}}/\sqrt t= \lambda/2$, so Lemma
\ref{lemma:1-ann} yields
\begin{eqnarray*}
I
&=&
\Big(\aver{\widehat{C}_1(\widetilde{B})}
|T_t(\bigchi_{\widetilde{B}}\,S_t f) |^{r}\,d\mu\Big)^{\frac1r}
\lesssim
\dec{{\lambda}}^{\theta_2}\,
\, e^{-c\,\lambda^2}\,\Big(\aver{\widetilde{B}}|S_t f|^{q}\,d\mu\Big)^{\frac1q}
\\
&=&
\dec{\lambda}^{\theta_2}\, e^{-c\,\lambda^2}
\,\Big(\aver{\widetilde{B}}|S_t(
f\bigchi_{\widetilde{B}})|^{q}\,d\mu\Big)^{\frac1q}
\lesssim
\dec{\lambda}^{\theta_2}\,e^{-c\,\lambda^2}\,
\dec{{\lambda}}^{\gamma_2}
\Big(\aver{\widetilde{B}} |f|^{p}\,d\mu\Big)^{\frac1p}
\\
&\lesssim&
\dec{\lambda}^{\theta_2+\gamma_2} \,e^{-c\,\lambda^2}\,
\Big(\aver{B} |f|^{p}\,d\mu\Big)^{\frac1p}.
\end{eqnarray*}
On the other hand,
\begin{eqnarray*}
II
&\lesssim&
\Big(\aver{2^{j+2}\,B}|T_t(
\bigchi_{2^{j+2}\,B\setminus 2^{j-1}\,B}\,S_t f)|^{r}\,d\mu
\Big)^{\frac1r}
+
\Big(\aver{2^{j+1}\,B}|T_t(
\bigchi_{(2^{j+2}\,B)^c}\,S_t f)|^{r}\,d\mu
\Big)^{\frac1r}
\\
&=&
II_1+II_2.
\end{eqnarray*}
For $II_1$ we use Lemma \ref{lemma:Bc}:
\begin{eqnarray*}
II_1
&=&
\Big(\aver{2^{j+2}\,B}|T_t(
\bigchi_{2^{j+2}\,B}\,(\bigchi_{(2^{j-1}\,B)^c}\,S_t f))|^{r}\,d\mu
\Big)^{\frac1r}
\\
&\lesssim&
\dec{\frac{r_{2^{j+2}\,B}}{\sqrt{t}}}^{\theta_2}\,
\,\Big(\aver{2^{j+2}\,B}\bigchi_{(2^{j-1}\,B)^c}\,|S_t
f|^{q}\,d\mu\Big)^{\frac1q}
\\
&\lesssim&
\dec{\lambda}^{\theta_2}\, 2^{(j-1)\gamma_1}\,
\dec{\frac{2^{j-1}\,r_B}{\sqrt{t}}}^{\max\{\gamma_1+D/q,\gamma_2,\}}\,
\expt{-\frac{c\,4^{j-1}\,r_B^2}{t}}\,
\Big(\aver{B} |f|^{p}\,d\mu\Big)^{\frac1p}
\\
&\lesssim&
2^{j\,\gamma_1}\, \dec{\lambda}^{\theta_2+
\max\{\gamma_1+D/q,\gamma_2,\}}\, e^{-c\,\lambda^2}\,
\Big(\aver{B} |f|^{p}\,d\mu\Big)^{\frac1p}.
\end{eqnarray*}
On the other hand for $II_2$ we set $\widetilde{B}=2^{j+1}\,B$ and so
its radius is $r_{\widetilde{B}}=2^{j+1}\,r$. Thus,
$(2^{j+2}\,B)^c=(2\,\widetilde{B})^c$ and $ r_{\widetilde{B}}/\sqrt t= 2\lambda$, so by Lemma \ref{lemma:Bc}
we have
\begin{eqnarray*}
\lefteqn{II_2
=
\Big(\aver{\widetilde{B}}
|T_t(\bigchi_{(2\,\widetilde{B})^c}\,S_t
f)|^{r}\,d\mu\Big)^{\frac1r}}
\\
&\lesssim&
2^{\theta_1}\,
\dec{\frac{2\,r_{\widetilde{B}}}{\sqrt{t}}}^{\max\{\theta_1,\theta_2\}}\,
e^{-c\,\lambda^2}
\,\Big(\aver{(2\,\widetilde{B})^c}|S_t f|^{q}\,d\mu\Big)^{\frac1q}
\\
&\lesssim&
\dec{\lambda}^{\max\{\theta_1,\theta_2\}}\, e^{-c\,\lambda^2}
\,\Big(\aver{(2^{j+2}\,B)^c}|S_t f|^{q}\,d\mu\Big)^{\frac1q}
\\
&\lesssim&
\dec{\lambda}^{\max\{\theta_1,\theta_2\}}\,e^{-c\,\lambda^2}\,
2^{(j+2)\,\gamma_1} \,
\dec{\frac{2^{j+2}\,r_B}{\sqrt{t}}}^{\max\{\gamma_1+D/q,\gamma_2\}}\,
\expt{-\frac{c\,4^{j+2}\,r_B^2}{t}}\,
\Big(\aver{\widetilde{B}} |f|^{p}\,d\mu\Big)^{\frac1p}
\\
&\lesssim&
2^{j\,\gamma_1}\,\dec{\lambda}^{\max\{\theta_1,\theta_2\}+\max\{\gamma_1+D/q,\gamma_2\}}
\,e^{-c\,\lambda^2}\,
\Big(\aver{B} |f|^{p}\,d\mu\Big)^{\frac1p}.
\end{eqnarray*}
Collecting the bounds for $I$, $II_1$ and $II_2$ we
obtain the desired estimate.
\end{proof}

\begin{proof}[Proof of Lemma \ref{lemma:sum}]
If $s\ge 1$, since $s^\beta\,e^{-c\,4^k\,s^2}\lesssim
e^{-c'\,s^2}\cdot e^{-c''\,4^k}$ for some $c''>0$, we have
$$
\sum_{k=0}^\infty
2^{k\,\alpha}\,\dec{2^k\,s}^{\beta}\,e^{-c\,4^k\,s^2}
\lesssim
e^{-c'\,s^2}\sum_{k=0}^\infty 2^{k\,(\alpha+\beta)}\,e^{-c''\,4^k}
\lesssim
e^{-c'\,s^2}
\lesssim
\dec{s}^{\max\{\alpha,\beta\}}\,e^{-c'\,s^2}.
$$
If $0<s<1$ then there is $k_0\in \NN$ such that $2^{-k_0}\le
s<2^{-k_0+1}$. We obtain
\begin{eqnarray*}
\lefteqn{\hskip-0.3cm
\sum_{k=0}^\infty
2^{k\,\alpha}\,\dec{2^k\,s}^{\beta}\,e^{-c\,4^k\,s^2}
\lesssim
 \sum_{k=0}^\infty
2^{k\,\alpha}\,\dec{2^{k-k_0}}^{\beta}\,e^{-c\,4^{k-k_0}}}
\\
&=&
\sum_{k=0}^{k_0}
2^{k\,\alpha}\,\dec{2^{k-k_0}}^{\beta}\,e^{-c\,4^{k-k_0}} +
\sum_{k=k_0+1}^\infty
2^{k\,\alpha}\,\dec{2^{k-k_0}}^{\beta}\,e^{-c\,4^{k-k_0}}
=
I+II.
\end{eqnarray*}
Then, as $\alpha\neq \beta$,
$$
I
\le
\sum_{k=0}^{k_0} 2^{k\,\alpha}\,2^{-(k-k_0)\,\beta}
\lesssim
2^{k_0\,\max\{\alpha,\beta\}}
\lesssim
s^{-\max\{\alpha,\beta\}}.
$$
On the other hand, as $\alpha+\beta>0$
$$
II
\le
\sum_{k=k_0+1}^\infty
2^{k\,\alpha}\,2^{(k-k_0)\,\beta}\,e^{-c\,4^{k-k_0}}
\lesssim
2^{k_0\,\alpha} \sum_{k=1}^\infty 2^{k\,(\alpha+\beta)}\,e^{-c\,4^k}
\lesssim
2^{k_0\,\alpha}
\lesssim
s^{-\max\{\alpha,\beta\}}.
$$
Thus,
$$
\sum_{k=0}^\infty
2^{k\,\alpha}\,\dec{2^k\,s}^{\beta}\,e^{-c\,4^k\,s^2}
\le
\,s^{-\max\{\alpha,\beta\}}
\lesssim
\dec{s}^{\max\{\alpha,\beta\}}\,e^{-c'\,s^2}.
$$
\end{proof}

\begin{proof}[Proof of Lemma \ref{lemma:1-ann}]
By Lemma \ref{lemma:geom:X}, given $B$ we can construct a sequence
$\{x_k\}_{k=1}^K\subset B$ with $K\le N^3$ such that $d(x_k, x_j)>
r_B/8$ for $j\ne k$ and with the property that for all $x\in B$ we
have some $k$ for which $d(x,x_k)\le r_{B}/8$ (this means that we
cannot pick more $x_k$'s). Write $B_k=B(x_k,r_{B}/4)$ and note that
$B\subset
\bigcup_{k=1}^{K} B_k$. Besides, $\widehat{C}_1(B)\subset 2^5\,B_k\setminus
2^2\,B_k=C_{2}(B_{k})\cup C_{3}(B_{k}) \cup C_{4}(B_{k})$ for each $k$. Let $f$ be supported on  $\widehat{C}_1(B)$.
Then, for each $k$, $f= \sum_{{j=2}}^4 f_{j,k}$ where
$f_{j,k}=f\,\bigchi_{C_{j}(B_{k})}$. As $\supp
f_{j,k}\subset C_{j}(B_k)$,
\begin{eqnarray*}
\lefteqn{\hskip-1.3cm
\Big(\aver{B}|T_t f|^{q}\,d\mu\Big)^{\frac1q}
\le
\sum_{k=1}^{K} \left(\frac{\mu(B_k)}{\mu(B)}\right)^{\frac1q}
\Big(\aver{B_k}|T_t f
|^{q}\,d\mu\Big)^{\frac1q}
\lesssim
\sum_{k=1}^{K}\sum_{j=2}^4
\Big(\aver{B_k}|T_t f_{j,k}
|^{q}\,d\mu\Big)^{\frac1q}}
\\
&\lesssim&
\sum_{k=1}^K\sum_{j=2}^4 2^{j\,\theta_1}\,
\dec{\frac{2^j\,r(B_k)}{\sqrt{t}}}^{\theta_2} \,
\expt{-\frac{c\,4^j\,r(B_k)^2}{t}} \,
\Big(\aver{C_j(B_k)}|f_{j,k}|^{p}\,d\mu\Big)^{\frac1p}
\\
&\lesssim&
\dec{\frac{2\,r}{\sqrt{t}}}^{\theta_2} \, \expt{-\frac{c\,4\,r^2}{t}}
\,
\Big(\aver{\widehat{C}_1(B)}|f|^{p}\,d\mu\Big)^{\frac1p}.
\end{eqnarray*}
where we have used that $\mu(2^{j+1}B_{k})\approx \mu(4B)\approx \mu(B)$ for $j=2,
3, 4$ and $1\le k \le K$.

On the other hand if $\supp f\subset B$ we have that $f=\sum_{k=1}^N
f_k$ where $f_k=f\,\bigchi_{E_k}$ with $E_k\subset B_k$ and the sets
$E_k$ are pairwise disjoint (for instance, we can take $E_{1}=B_{1}$,
$E_{2}=B_{2}\setminus E_{1}$, \dots). Then, as $\supp f_k\subset B_k$
\begin{eqnarray*}
\lefteqn{\hskip-2cm
\Big(\aver{\widehat{C}_1(B)}|T_t f|^{q}\,d\mu\Big)^{\frac 1 q}
\le
\sum_{k=1}^K\sum_{j=2}^4
\left(\frac{\mu(2^{j+1}\,B_k)}{\mu(4\,B)}\right)^{\frac 1 q}
\Big(\aver{C_j(B_k)}
|T_t f_k|^{q}\,d\mu\Big)^{\frac 1 q}}
\\
&\lesssim&
\sum_{k=1}^K\sum_{j=2}^4 2^{j\,\theta_1}\,
\dec{\frac{2^j\,r(B_k)}{\sqrt{t}}}^{\theta_2}
\,\expt{-\frac{c\,4^j\,r(B_k)^2}{t}} \,
\Big(\aver{B_k}|f_k|^{p}\,d\mu\Big)^{\frac 1 p }
\\
&\lesssim&
\dec{\frac{2\,r}{\sqrt{t}}}^{\theta_2} \, \expt{-\frac{c\,4\,r^2}{t}}
\,
\Big(\aver{B}|f|^{p}\,d\mu\Big)^{\frac 1 p }.
\end{eqnarray*}
\end{proof}

\begin{proof}[Proof of Lemma \ref{lemma:Bc}]
Suppose first that $j\ge 2$. Let us write $\lambda=2^j\,r/\sqrt{t}$.
Using that $T_t\in \off{p}{q}$ we have
\begin{eqnarray*}
\lefteqn{\hskip-1cm
\Big(\aver{B}|T_t( \bigchi_{(2^j\,B)^c}\, f) |^{q}\,d\mu\Big)^{\frac 1 q}
\le
\sum_{k\ge j}
\Big(\aver{B}|T_t( \bigchi_{C_k(B)}\, f) |^{q}\,d\mu\Big)^{\frac 1 q}}
\\
&\lesssim&
\sum_{k\ge j}2^{k\,\theta_1}
\dec{\frac{2^k\,r}{\sqrt{t}}}^{\theta_2}\,
\expt{-\frac{c\,4^k\,r^2}{t}} \,
\Big(\aver{C_k(B)}|f|^{p}\,d\mu\Big)^{\frac 1 p }
\\
&\le&
2^{j\,\theta_1}\,\sum_{k=0}^\infty 2^{k\,\theta_1}
\dec{2^k\,\lambda}^{\theta_2}\, e^{-c\,4^k\,\lambda^2} \,
\Big(\aver{(2^j\,B)^c}|f|^{p}\,d\mu\Big)^{\frac 1 p }
\\
&\lesssim&
2^{j\,\theta_1}\,\dec{\lambda}^{\max\{\theta_1,\theta_2\}}\,
e^{-c'\,\lambda^2} \,
\Big(\aver{(2^j\,B)^c}|f|^{p}\,d\mu\Big)^{\frac 1 p },
\end{eqnarray*}
where we have used Lemma \ref{lemma:sum} (the power
$\max\{\theta_1,\theta_2\}$ is correct whenever
$\theta_1\neq\theta_2$, see Remark \ref{remark:powers-diff}
otherwise).

When $j=1$, the argument is exactly the same but for the term $k=j=1$
on which we use Lemma \ref{lemma:1-ann} in place of Definition
\ref{defi:off-d:weights}.

On the other hand, assume that $j\ge 2$ (and the other case is done
as just explained). Then, since $\mu$ is doubling, we have
\begin{eqnarray*}
\lefteqn{\hskip-1.3cm
\Big(\aver{(2^j\,B)^c}|T_t( \bigchi_{B}\, f) |^{q}\,d\mu\Big)^{\frac 1 q}
\le
\sum_{k\ge j}
\left(\frac{\mu(2^{k+1}\,B)}{\mu(2^{j+1}\,B)}\right)^{\frac 1 q}
\Big(\aver{C_k(B)}|T_t( \bigchi_{B}\, f) |^{q}\,d\mu\Big)^{\frac 1 q}}
\\
&\lesssim&
\sum_{k\ge j}2^{(k-j)\,D/q}\, 2^{k\,\theta_1}\,
\dec{\frac{2^k\,r}{\sqrt{t}}}^{\theta_2}\,
\expt{-\frac{c\,4^k\,r^2}{t}} \,
\Big(\aver{B}|f|^{p}\,d\mu\Big)^{\frac 1 p }
\\
&=&
2^{j\,\theta_1}\,\sum_{k=0}^\infty 2^{k\,(\theta_1+D/q)}\,
\dec{2^k\,\lambda}^{\theta_2}\, e^{-c\,4^k\,\lambda^2} \,
\Big(\aver{B}|f|^{p}\,d\mu\Big)^{\frac 1 p }
\\
&\lesssim&
2^{j\,\theta_1}\,\dec{2^k\,\lambda}^{\max\{\theta_1+D/q,\theta_2\}}\,
e^{-c'\,4^j\,\lambda^2} \,
\Big(\aver{(2^j\,B)^c}|f|^{p}\,d\mu\Big)^{\frac 1 p },
\end{eqnarray*}
where, as before, we have used Lemma \ref{lemma:sum} (here, one needs
$\theta_1+D/q\neq \theta_2$, see Remark \ref{remark:powers-diff}
otherwise).
\end{proof}

\subsection{Proof of Proposition \ref{prop:off-unw-w}}
\label{sec:proof:prop:off-unw-w}

We fix $w\in A_{\frac{p}{p_0}}\cap  RH_{(\frac {q_{0}}q)'}$ and by
Proposition \ref{prop:weights}, $(iii)$ and $(iv)$, there exist $p_1$, $q_1$ with
$p_0<p_1<p\le q<q_1<q_0$ such that $w\in A_{\frac{p}{p_1}}\cap
RH_{(\frac {q_{1}}q)'}$.  It is well-known that
\begin{eqnarray}
w \in A_{\frac{p}{p_1}}
&\Longleftrightarrow&
\Big (\aver{B} g^{p_{1}} \, d\mu \Big)^{\frac 1 {p _{1}}}
\lesssim
 \Big(\aver{B} g^{p}\,  dw\Big)^{\frac 1 p } \label{Ap-g}
\\[0.2cm]
w \in RH_{(\frac {q_{1}}q)'}
&\Longleftrightarrow&
\Big(\aver{B} g^q\, dw\Big)^{\frac 1 q}
\lesssim
  \Big(\aver{B} g^{q_{1}}\,  d\mu \Big)^{\frac 1 {q_{1}}}, \label{RH-g}
\end{eqnarray}
where in the right hand sides $g$ runs over the set of non-negative
measurable functions and $B$ runs over the set of balls. Thus, using
\eqref{RH-g}, $T_t\in\off{p_1}{q_1}$ and \eqref{Ap-g} we have
\begin{align*}
\Big(\aver{B} |T_t( \bigchi_B\, f) |^{q}\,dw\Big)^{\frac 1 q}
&\lesssim
\Big(\aver{B} |T_t( \bigchi_B\, f) |^{q_1}\,d\mu\Big)^{\frac 1 {q_1}}
\lesssim
\dec{\frac{r}{\sqrt{t}}}^{\theta_2} \,\Big(\aver{B}
|f|^{p_1}\,d\mu\Big)^{\frac 1 {p _1}}
\\
&\lesssim
\dec{\frac{r}{\sqrt{t}}}^{\theta_2} \,\Big(\aver{B} |f|^{p}\,d
w\Big)^{\frac 1 p }.
\end{align*}
This shows \eqref{w:off:B-B}. The  same can be done to derive
\eqref{w:off:C-B} and \eqref{w:off:B-C} and this completes the proof.

\subsection{Proof of Propositions \ref{prop:full-off} and \ref{prop:full-off:2D}}\label{sec:proof:prop:full-off}

\begin{proof}[Proof of Proposition  \ref{prop:full-off}] We prove $(a)$.  That $T_{t}\in \full{p}{q}$ implies $T_{t}\in
\off{p}{q}$ is easy by specializing \eqref{eq:offLpLq-w} to balls and
annuli and using
$$
{\mu(B)^{\frac 1 p -\frac 1 q} } \lesssim r^{\frac n p - \frac n q}.
$$
from  the  polynomial upper bound of the volume.

We turn to $(b)$.  Assume $q<\infty$. The argument mimics that of
Theorem \ref{theor:uniform-compo}, part $(a)$. Let $E,F$ be two
closed sets and $t>0$. Let $f$ be supported in $E$.  We first assume
that $ t< (d(E,F)/16)^2$. Pick a collection of balls  $B_k=B(x_k,r)$
as in Lemma \ref{lemma:covering-X} with $r=d(E,F)/16$.  Observe that
if $x\in F$ and $y\in E$ then $d(x,y)\ge d(E,F)= 16\,  r$. Hence, if
$x\in B_{k}$ then $y\notin 4B_{k}$, so $y\in C_{j}(B_{k})$ for some
$j\ge 2$. In what follows the summation in $k$ is restricted to those
balls $B_{k}$ so that $F\cap B_k\neq\emptyset$. Using that $\supp f
\subset E$ and \eqref{w:off:C-B}, we have
\begin{eqnarray*}
\lefteqn{\|T_t f\|_{L^q(F,\mu)}^q
\le
\sum_{k} \int_{B_k}|T_t f|^q\, d\mu
\le
\sum_{k} \bigg( \sum_{j=2}^\infty
\Big(\int_{B_k}|T_t (\bigchi_{C_j(B_k)}\, f)|^q\, d\mu\Big)^{\frac 1 q}
\bigg)^{q}}
\\
&\lesssim&
\sum_{k} \bigg( \sum_{j=2}^\infty
2^{j\,\theta_1}\,\dec{\frac {2^j r }{\sqrt t}}^{\theta_2}\, e^{-\frac{c\,4^j\, r^2}{t}}\,
\frac{\mu(B_k)^{\frac 1 q}}{\mu(2^{j+1}\,B_k)^{\frac 1 p }}
\Big(\int_{C_j(B_k)}|f|^p\, d\mu\Big)^{\frac 1 p }
\bigg)^{q}.
\end{eqnarray*}
Next,  by the  polynomial lower bound of the volume, $p\le q$ and $r > \sqrt t$,  we have
$$
\frac{\mu(B_k)^{\frac 1 q}}{\mu(2^{j+1}\,B_k)^{\frac 1 p } }\le
\frac{\mu(B_k)^{\frac 1 q}}{\mu(B_k)^{\frac 1 p }}
\lesssim r^{\frac n q - \frac n p}
\le
t^{-\frac12\,(\frac{n}{p} - \frac{n}{q})} .
$$
Also,  note that
$$
2^{j\,\theta_1}\,\dec{\frac {2^j r }{\sqrt t}}^{\theta_2}\, e^{-\frac{c\,4^j\, r^2}{t}} \lesssim e^{-c'\,4^j}\, e^{-\frac{c''\, r^2}{t}}
$$
for some $c',c''>0$. Thus,
 \begin{eqnarray*}
\|T_t f\|_{L^q(F,\mu)}^q
&\lesssim&
\bigg( t^{-\frac12\,(\frac{n}{p} - \frac{n}{q})} \,  e^{-\frac{c'' r^2}{t}}\bigg)^q
\sum_{k} \bigg( \sum_{j=2}^\infty e^{-c'\,4^j}\,
\Big(\int_{C_j(B_k)}|f|^p\, d\mu\Big)^{\frac 1 p }\bigg)^q
\\
&\lesssim&
\bigg( t^{-\frac12\,(\frac{n}{p} - \frac{n}{q})} \,  e^{-\frac{c'' r^2}{t}}\bigg)^q
\sum_{k} \bigg( \sum_{j=2}^\infty e^{-c'\,4^j}\,
\int_{C_j(B_k)}|f|^p\, d\mu\bigg)^{\frac q p}
\\
&\le&
\bigg( t^{-\frac12\,(\frac{n}{p} - \frac{n}{q})} \,  e^{-\frac{c'' r^2}{t}}\bigg)^q
\bigg(\sum_{k}  \sum_{j=2}^\infty e^{-c'\,4^j}\,
\int_{C_j(B_k)}|f|^p\, d\mu\bigg)^{\frac q p}
\\
&=&
\bigg( t^{-\frac12\,(\frac{n}{p} - \frac{n}{q})} \,  e^{-\frac{c'' r^2}{t}}\bigg)^q
\bigg(
\int_{E}\sum_{j=2}^\infty \sum_{k}   e^{-c'\,4^j}\,
\bigchi_{C_j(B_k)}\, |f|^p\, d\mu\bigg)^{\frac q p}
\end{eqnarray*}
where we used the fact that $\frac q p \ge 1$.
We conclude as in Section \ref{sec:proof:theor:uniform-compo-a}
using
$$
\sum_{j=2}^\infty\sum_{k}   e^{-c'\,4^j}\, \bigchi_{C_j(B_k)} \le
\sum_{j=2}^\infty N^{j+3} e^{-c'\, 4^j} \le C <\infty.
$$

In the case where $t\ge  (d(E,F)/16)^2$, then we argue as before with
$r=\sqrt t$. This time, we have to incorporate the terms with $j=1$
and use also  \eqref{w:off:B-B} with $4B_{k}$.  We obtain
$$
\|T_t f\|_{L^q(F,\mu)} \lesssim  t^{-\frac12\,(\frac{n}{p} -
\frac{n}{q})}\,  \|f\|_{L^p(E,\mu)}.
$$

This proves the result when $q<\infty$. The modification for
$q=\infty$ is left to the reader.
\end{proof}

\begin{proof}[Proof of Proposition  \ref{prop:full-off:2D}]

Assume first that $T_{t}$ satisfies $L^p(\mu)-L^q(\mu)$ full
off-diagonal estimates of type $\varphi$.   Let us see for example,
how to obtain \eqref{w:off:B-B} for $T_{t}\in \off{p}{q}$, the other
estimates being similar.  We  specialize \eqref{eq:offLpLq-w} with
$t^{-\theta}$ replaced by $\varphi\big({\sqrt t}\big)^{\frac 1 q -\frac 1 p}$
to $E=F=B$. Then, using $\mu(B)\sim \varphi(r)$, we obtain
$$
\Big(\aver{B} |T_t f|^q\, d\mu\Big)^{\frac 1 q}
\lesssim
\left( \frac {\varphi(r)}{\varphi\big(\sqrt t\big)}\right)^{\frac 1 p -\frac 1 q} \,
\Big(\aver{B} |f|^p\, d\mu\Big)^{\frac 1 p}
$$
and we observe that since $\varphi$ is non decreasing and $\mu$ is doubling,
$$
\frac {\varphi(r)}{\varphi\big(\sqrt t\big)}
\lesssim
\max\Big\{ \Big(\frac r {\sqrt t}\Big)^D, 1\Big\}
\lesssim
\dec{ \frac r {\sqrt t}}^{ D}
$$
where $D$ is the doubling order of $\mu$.

We turn to the converse. The argument is the same as the one above.
The only change is in the inequality when $t <(d(E,F)/16)^2=r^2$, and it reads
$$
\frac{\mu(B_k)^{\frac 1 q}}{\mu(2^{j+1}\,B_k)^{\frac 1 p } }\le
\frac{\mu(B_k)^{\frac 1 q}}{\mu(B_k)^{\frac 1 p }} \sim
\varphi(r)^{\frac 1 q - \frac 1 p}
\le
\varphi\big(\sqrt t\big)^{\frac{1}{q} - \frac{1}{p}}
$$
as $\varphi$ is non-decreasing  $\sqrt t \le r$ and $p\le q$. Further details are left to the reader.
\end{proof}

\subsection{Proof of Proposition \ref{prop:full-strong}}\label{sec:proof:prop:full-strong}

 Fix a ball $B$. Set $r=r_{B}$ and $\lambda=r/\sqrt{t}$. By $(iii)$
and $(iv)$ in Proposition \ref{prop:weights}, one can select $p_1,
q_{1}$ with $p_0<p_1<p\le q<q_1<q_0$ and $w \in A_{\frac p {p_{1}}}
\cap RH_{(\frac {q_{1}}q)'}$.

Using $(vii)$ and $(viii)$ in Proposition \ref{prop:weights},  we
have that $w^{(\frac  {q_{1}}q)'} \in A_{\alpha}$ with $\alpha
=1 + (\frac p {p_{1}}-1) (\frac  {q_{1}}q)'$ and
\begin{align*}
\bigg( \int _{(2B)^c}  w^{(\frac  {q_{1}}q)'}
\left( \frac {|x-x_{B}|}{r}\right)^{-n\alpha} dx\bigg)^{1/(\frac
{q_{1}}q)'}
&
\lesssim
\left( \int _{B}  w^{(\frac  {q_{1}}q)'}(x)
dx\right)^{1/(\frac  {q_{1}}q)'}
\lesssim
\frac{w(B)}{|B|^{\frac q {q_{1}}}},
\end{align*}
where, in the last estimate we have used that $w \in RH_{(\frac
{q_{1}}q)'}$. Let $a>0$ be such that $n\alpha = a (\frac {q_{1}}q)'
$, then by H\"older's inequality and the above inequality
$$
\Big(\aver{(2B)^c} |T_{t}(\bigchi_{B}\, f)|^q \, dw \Big)^{\frac1q} \lesssim
\Big(\aver{(2B)^c} |T_{t}( \bigchi_{B}\, f)|^{q_{1}}\, \bigg( \frac
{|x-x_{B}|}{r}\bigg)^{a \frac {q_{1}} q}\, dx
\Big)^{\frac{1}{q_{1}}}.
$$
Now decompose $(2B)^c$ as the union of the rings $C_{j}(B)$ for $j\ge
1$ where $C_{1}(B)$ is here $4B\setminus 2B$. On each $C_{j}(B)$ we
can use the $L^{p_{1}}(dx)-L^{q_{1}}(dx)$ full off-diagonal estimates
so that the right hand term is bounded by
$$
\Big(\sum_{j\ge 1} 2^{j\,a\, \frac {q_{1}} q} e^{-c\,4^j\,\lambda^2}
\Big)^{\frac1{q_{1}}}\,
\lambda^{\frac{n}{p_{1}} - \frac{n}{q_{1}}}\,
\Big(\aver{B} |f|^{p_{1}} \, dx \Big)^{\frac1{p_{1}}}
\lesssim
\lambda^{\frac{n}{p_{1}} - \frac{n}{q_{1}}-\frac {a}{q}}\,
e^{-c\,\lambda^2} \,
\Big(\aver{B}
|f|^{p} \, dw \Big)^{\frac1p}
$$
where we have used that $w\in A_{\frac p {p_{1}}}$. Hence we have
obtained \eqref{eq:full-strong:B-2Bc} and the total power of
$\lambda$ is
$$
\beta= \frac n {p_{1}} - \frac n{q_{1}} - \frac a {q}= \frac{
n}{p_{1}}\left( 1 - \frac p q \right).
$$

Let us prove the estimate \eqref{eq:full-strong:2Bc-B}. We pick
$p_{1},q_{1}$ as before and take $a$ so that $a (\frac p
{p_{1}})'=n\alpha'$ where $\alpha =1 + (\frac p {p_{1}}-1) (\frac
{q_{1}}q)'$. Writing $f_{j}=\bigchi_{C_{j}(B)} f$, using that $w\in
RH_{\left(\frac{q_1}{q}\right)'}$ and by the
$L^{p_{1}}(dx)-L^{q_{1}}(dx)$ full off-diagonal estimates for $T_{t}$
\begin{align*}
\lefteqn{\hskip-.1cm
\Big(\aver{B} |T_{t}(\bigchi_{(2B)^c}\, f)|^q\, dw \Big)^{\frac1q}
\lesssim
\Big(\aver{B} |T_{t}(\bigchi_{(2B)^c}\, f)|^{q_{1}}\,  dx
\Big)^{\frac1{q_{1}}}
\le
\sum_{j\ge 1} \Big(\aver{B} |T_{t}f_{j}|^{q_{1}}\,  dx
\Big)^{\frac1{q_{1}}}}
\\
&
\lesssim
\frac1{|B|^{\frac1{p_1}}}\, \sum_{j\ge 1} \lambda^{\frac{n}{p_{1}} -
\frac{n}{q_{1}}}\, e^{-c\, 4^j\,\lambda^2}\,
\Big(\int_{C_{j}(B)} |f|^{p_{1}}\,  dx \Big)^{\frac1{p_{1}}}.
\\
&
\lesssim
\frac1{|B|^{\frac1{p_1}}}\sum_{j\ge 1} \lambda^{\frac{n}{p_{1}} -
\frac{n}{q_{1}}}\, e^{-c\, 4^j\,\lambda^2}\, 2^{j\,\frac{a}{p_{1}}}
\Big(\int_{C_{j}(B)} \Big( \frac
{|x-x_{B}|}{r}\Big)^{-a} |f|^{p_{1}}\,  dx \Big)^{\frac1{p_{1}}}
\\
&
\le
\frac1{|B|^{\frac1{p_1}}}
\Big(
\sum_{j\ge 1}
\Big[
\lambda^{\frac{n}{p_{1}} - \frac{n}{q_{1}}}\,
e^{-c\,4^j\,\lambda^2}\, 2^{j\,\frac{a}{p_{1}}}
\Big]^{p_1'}
\Big)^{\frac1{p_1'}}
\,
\Big(\sum_{j\ge 1}\int_{C_{j}(B)} \Big( \frac
{|x-x_{B}|}{r}\Big)^{-a}|f|^{p_{1}} \,  dx \Big)^{\frac1{p_{1}}}
\\
&
\lesssim
\frac{\lambda^{\gamma}\, e^{-c\,\lambda^2}}{|B|^{\frac1{p_1}}}\,
\Big(\int_{(2B)^c} \Big( \frac {|x-x_{B}|}{r}\Big)^{-a} |f|^{p_{1}} \,
dx
\Big)^{\frac1{p_{1}}},
\end{align*}
where
$$\gamma = \frac n {p_{1}} - \frac  n {q_{1}} - \frac a {p_{1}}=
\frac n {q_{1}}
\left( \frac q p - 1 \right).
$$
Using H\"older's inequality with $\frac p {p_{1}}>1$ it follows that
\begin{align*}
\lefteqn{
\Big(\aver{B} |T_{t}(\bigchi_{(2B)^c}\, f)|^q\, dw \Big)^{\frac1q}}
\\
&\lesssim
\frac{\lambda^{\gamma}\, e^{-c\,\lambda^2}}{|B|^{\frac1{p_1}}}\,
\Big(\int_{(2B)^c} |f|^{p}w \,  dx \Big)^{\frac1p}\,
\Big(\int_{(2B)^c} w(x)^{1-(\frac p {p_{1}})'}
\Big(\frac {|x-x_{B}|}{r}\Big)^{-a(\frac p {p_{1}})'} \,  dx
\Big)^{\frac{1}{p_1\,(\frac p {p_{1}})'}}
\\
&
\lesssim
\lambda^{\gamma}\, e^{-c\,\lambda^2}\,
\Big(\aver{(2B)^c} |f|^{p} \,dw\Big)^{\frac1p}\,
\frac{w(B)^\frac1p}{|B|^\frac1{p_1}}
\Big(\int_{(2B)^c} w^{1-(\frac p {p_{1}})'}(x)
\Big(\frac {|x-x_{B}|}{r}\Big)^{-n\,\alpha'} \,  dx
\Big)^{\frac{1}{p_1\,(\frac p {p_{1}})'}}
\\
&
\lesssim
\lambda^{\gamma}\, e^{-c\,\lambda^2}\,\Big(\aver{(2B)^c} |f|^{p}
\,dw\Big)^{\frac1p},
\end{align*}
where the latter inequality is obtained as follows: since $w^{(\frac
{q_{1}}q)'} \in A_{\alpha}$, one has that $w^{(\frac
{q_{1}}q)'(1-\alpha')}= w^{1- (\frac p {p_{1}})'}\in A_{\alpha'}$.
Thus $(viii)$ in Proposition \ref{prop:weights} and $w \in A_{\frac p
{p_{1}}}$ imply
\begin{align*}
\lefteqn{\hskip-3cm
\bigg(\int _{(2B)^c} w^{1-(\frac  p{p_{1}})'}(x)
\Big( \frac{|x-x_{B}|}{r}\Big)^{-n\alpha'} dx
\bigg)^{\frac{1}{p_1\,(\frac p
{p_{1}})'}}
\lesssim
\Big(\int_{B} w^{1-(\frac  p{p_{1}})'}(x) \,
dx\Big)^{\frac{1}{p_1\,(\frac p {p_{1}})'}}}
\\
&
\lesssim
|B|^{\frac{1}{p_1\,(\frac p {p_{1}})'}}\,
\left(\frac{|B|}{w(B)}\right)^{\frac{(\frac{p}{p_{1}})'-1}{p_1\,(\frac{p}{p_{1}})'}}
=
\frac {|B|^\frac1{p_1}}{w(B)^{\frac1p}}.
\end{align*}

\appendix\section{Muckenhoupt weights}\label{App:Weights}

A weight $w$ is a non-negative locally integrable function. We say
that $w\in A_p$, $1<p<\infty$, if there exists a constant $C$ such
that for every ball $B\subset\X$,
$$
\Big(\aver{B} w\,d\mu\Big)\,
\Big(\aver{B} w^{1-p'}\,d\mu\Big)^{p-1}\le C.
$$
For $p=1$, we say that $w\in A_1$ if there is a constant $C$ such
that for every ball $B\subset \X$
$$
\aver{B} w\,d\mu
\le
C\, w(y),
\qquad \mbox{for a.e. }y\in B,
$$
or, equivalently, $M w\le C\,w$ for a.e.. The reverse
H\"older classes are defined in the following way: $w\in RH_{q}$, $1<
q<\infty$, if
$$
\Big(\aver{B} w^q\,d\mu\Big)^{\frac1q}
\le C\, \aver{B} w\,d\mu
$$
for every ball $B$. The endpoint $q=\infty$ is given by the
condition: $w\in RH_{\infty}$ whenever, for any ball $B$,
$$
w(y)\le C\, \aver{B} w\,d\mu,
\qquad \mbox{for a.e. }y\in B.
$$
Notice that we have excluded the case $q=1$ since the class $RH_1$
consists of all the weights and that is the way $RH_1$ is understood.

Next, we present some of the properties of these classes.
\begin{prop}\label{prop:weights}\
\begin{enumerate}
\renewcommand{\theenumi}{\roman{enumi}}
\renewcommand{\labelenumi}{$(\theenumi)$}
\addtolength{\itemsep}{0.2cm}

\item $A_1\subset A_p\subset A_q$ for $1\le p\le q<\infty$.

\item $RH_{\infty}\subset RH_q\subset RH_p$ for $1<p\le q\le \infty$.

\item If $w\in A_p$, $1<p<\infty$, then there exists $1<q<p$ such
that $w\in A_q$.

\item If $w\in RH_q$, $1<q<\infty$, then there exists $q<p<\infty$ such
that $w\in RH_p$.

\item $\displaystyle A_\infty=\bigcup_{1\le p<\infty} A_p=\bigcup_{1<q\le
\infty} RH_q $

\item If $1<p<\infty$, $w\in A_p$ if and only if $w^{1-p'}\in
A_{p'}$.

\item If $1\le q\le \infty$ and $1\le s<\infty$, then $\displaystyle
w\in A_q \cap RH_s$ if and only if $ w^{s}\in A_{s\,(q-1)+1}$.

\item In the case of the Euclidean space $\RR^n$, if $w\in A_p$, $1<p<\infty$, there exists $C_w$ such that for
every ball $B=B(x_B,r_B)$,
$$
\int_{\RR^n\setminus 2\,B}
w(x)\,\left(\frac{|x-x_B|}{r_{B}}\right)^{-n\,p}\,dx
\le
C_w\,w(B).
$$

\end{enumerate}
\end{prop}

Properties $(i)$-$(vi)$ are standard. For $(vii)$ see \cite{JN} in
the Euclidean setting (and the same argument holds in  spaces of
homogeneous type \cite{ST}). The last property follows easily by
using the boundedness of $M$ on $L^{p}(w)$ applied to $f=\bigchi_B$.

\end{document}